\RequirePackage{fix-cm}
\documentclass{svjour3}                     %
\smartqed  %

\usepackage{tabularx}
\usepackage{amsmath}
\usepackage{amssymb}
\usepackage{dsfont}
\usepackage{rotating}
\usepackage{multirow}
\usepackage{epstopdf}
\usepackage{newtxmath}
\usepackage{mathtools}

\usepackage{multirow} %
\usepackage{longtable,tabu} %
\usepackage{hhline}
\usepackage[linesnumbered,ruled,vlined]{algorithm2e}

\usepackage{tikz}
\usetikzlibrary{math}
\usepackage{pgfplots}
\usepackage{siunitx}
\usepackage{nicefrac}

\usepackage{graphicx}
\graphicspath{{figures/}}

\usepackage{hyperref}
\hypersetup{
	colorlinks=true,
	linkcolor=black, %
	citecolor=black, %
	urlcolor=blue  } %

\usepackage{todonotes}

\DeclareMathOperator\erf{erf}

\newcommand{\R}{\mathbb{R}}
\newcommand{\Ogrande}{\mathcal{O}}

\renewcommand{\L}{\mathrm{L}}
\newcommand{\sign}{{\rm sign}}

\newcommand{\CFL}{\ensuremath{\text{CFL}}}
\newcommand{\SCOUT}{\ensuremath{\text{SCOUT}}}
\newcommand{\SL}{\ensuremath{\text{SL}}}
\newcommand{\SCOUTCN}{\ensuremath{\text{SCOUT-CN}}}
\newcommand{\SCOUTSLCN}{\ensuremath{\text{SCOUT-SLCN}}}
\newcommand{\dt}{\mathrm{\Delta}t}
\newcommand{\dx}{\mathrm{\Delta}x}
\newcommand{\tn}{\ensuremath{t^n}}
\newcommand{\tnn}{\ensuremath{t^{n+1}}}

\newcommand{\xipm}{\ensuremath{x_{i+1/2}}}
\newcommand{\ximm}{\ensuremath{x_{i-1/2}}}
\newcommand{\ipm}{\ensuremath{i+1/2}}
\newcommand{\imm}{\ensuremath{i-1/2}}
\newcommand{\pder}[2]{{\partial_{#1}#2}}
\newcommand{\PDER}[2]{\frac{\partial #2}{\partial #1}}

\newcommand{\intd}{\ensuremath{\mathrm{d}}}    %
\renewcommand{\vec}[1]{\ensuremath{\mathbf{#1}}}
\newcommand{\domainlag}{\Omega_L}
\newcommand{\rec}{\ensuremath{\mathcal{R}}}

\newfont{\numerikEleven}{ecrm1000}
\newfont{\numerikTen}{cmss10}
\newfont{\numerikNine}{cmss9}
\newfont{\numerikEight}{cmss8}

\begin{document}

\title{SCOUT: Semi-Lagrangian COnservative and Unconditionally sTable schemes for nonlinear advection-diffusion problems}

\titlerunning{Semi-Lagrangian conservative schemes for advection-diffusion problems}        %

\author{Silvia Preda \and Walter Boscheri \and Matteo Semplice \and Maurizio Tavelli}

\institute{S. Preda \at Department of Science and High Technology, University of Insubria, 22100 Como, Italy \\ \email{silvia.preda@uninsubria.it} 
	       \and	
	       W. Boscheri \at
	       Laboratoire de Mathématiques UMR 5127 CNRS, Universit{\'e} Savoie Mont Blanc, 73376 Le Bourget du Lac, France \\
	       Department of Mathematics and Computer Science, 44121 Ferrara, Italy \\	\email{walter.boscheri@univ-smb.fr}
   		   \and
   		   M. Semplice \at Department of Science and High Technology, University of Insubria, 22100 Como, Italy \\ \email{matteo.semplice@uninsubria.it} 
   		   \and
		   M. Tavelli \at
           Department of Engineering for Innovation Medicine, University of Verona, 37134 Verona, Italy  \\
		   \email{maurizio.tavelli@univr.it}
}

\date{Received: date / Accepted: date}

\maketitle

\begin{abstract}
In this work, we propose a new semi-Lagrangian ($\SL$) finite difference scheme for nonlinear advection-diffusion problems. To ensure conservation, which is fundamental for achieving physically consistent solutions, the governing equations are integrated over a space-time control volume constructed along the characteristic curves originating from each computational point. By applying Gauss theorem, all space-time surface integrals can be evaluated. For nonlinear problems, a nonlinear equation must be solved to find the foot of the characteristic, while this is not needed in linear cases. This formulation yields $\SL$ schemes that are fully conservative and unconditionally stable, as verified by numerical experiments with $\CFL$ numbers up to 100. Moreover, the diffusion terms are, for the first time, directly incorporated within a conservative semi-Lagrangian framework, leading to the development of a novel characteristic-based Crank-Nicolson discretization in which the diffusion contribution is implicitly evaluated at the foot of the characteristic. A broad set of benchmark tests demonstrates the accuracy, robustness, and strict conservation property of the proposed method, as well as its unconditional stability.
\subclass{MSC 65 \and MSC 68}
\end{abstract}

\section{Introduction} \label{sec:intro}

Advection and diffusion processes play a crucial role in describing how quantities like heat, pollutants, or nutrients are transported and spread within fluids, capturing the combined effects of large-scale flow (advection) and molecular mixing (diffusion). Indeed, solving advection-diffusion equations holds significant importance for society, as they underpin a wide range of applications which span from geophysical flows in oceans, rivers, and lakes, to problems in aerospace and mechanical engineering, atmospheric modeling for weather prediction, and even blood flow within the human cardiovascular system. The numerical solution of advection-diffusion problems therefore remains an active and significant area of research.

Semi-Lagrangian ($\SL$) methods have recently attracted considerable attention owing to their superior resolution and stability properties compared to more classical upwind finite difference and finite volume methods. In the $\SL$ framework, the advection term is reformulated in a Lagrangian setting and subsequently discretized. Specifically, $\SL$ schemes involve integrating material trajectories backward in time to determine the foot of the characteristic, where the numerical solution is evaluated through interpolation. Originally introduced in the context of numerical weather prediction \cite{SL_29,SL_30}, these methods have since been successfully extended to a broad range of applications, including environmental engineering for the simulation of free-surface flows in rivers and oceans \cite{SL_39,SL_40,ADERFSE}, plasma physics \cite{Vlasov2003,Vlasov2010,respaud2011,liu2025asymptotic}, kinetic theory \cite{Carrillo2007,SLBoltz12,SLkinI_Bos2021,boscarino2024conservative}, image processing \cite{Carlini2013,preda2025surface}, particle tracking in biology \cite{shotr2024}, and the numerical solution of Hamilton-Jacobi equations \cite{Falcone2002,carlini2025cweno}.

Semi-Lagrangian discretization techniques have been widely utilized in the context of finite difference schemes for modeling environmental applications \cite{Casulli1990,Casulli1999} as well as for achieving high-order interpolation accuracy in finite-volume formulations \cite{BoscheriDumbser,ader_fse}. In \cite{ELLAM2006}, a semi–Lagrangian approach for simulating solute transport and sorption processes in porous media is proposed. Staggered mesh configurations were explored in \cite{Bonaventura2000,rosatti2005semi} for the construction of $\SL$ algorithms. Additionally, discontinuous Galerkin (DG) formulations incorporating $\SL$ schemes for advection terms were presented in \cite{TumoloBonaventuraRestelli,TumoloBonaventura}. More recently, the semi-Lagrangian strategy has been applied to the incompressible Navier-Stokes equations in the vorticity-stream function formulation \cite{bonaventura2018vortexstream}, achieving unconditional stability for both advection and diffusion processes thanks to a semi-implicit time discretization. Parabolic problems have been treated in \cite{bonaventura2014siam,bonaventura2021advdiff} by means of a semi-Lagrangian schemes which also accounts for the diffusion effects at the foot of the characteristic by means of a Brownian motion approach.

A distinctive advantage of the semi-Lagrangian approach lies in its ability to employ large time steps without being restricted by a Courant-Friedrichs-Lewy ($\CFL$) type stability condition, thus providing an appealing alternative to conventional explicit upwind discretizations. This explains the remarkable growth in their development over the past decade, particularly for applications involving incompressible flows and linear transport problems or, more broadly, mathematical models where strict conservation is not a primary requirement. Nevertheless, the extension of $\SL$ methods to problems involving shock waves in advection equations remains nontrivial, as conservation must be rigorously preserved in order to obtain physically reliable solutions. Conservative formulations of semi-Lagrangian schemes addressing this issue can be found, for instance, in \cite{QuiShu2011,LentineEtAl2011,EulLagWENO2012,ConsSL2021,bergami2023}. The core idea proposed in \cite{QuiShu2011} is to integrate the numerical solution over a space-time control volume, which is obtained for each grid point by a backward trajectory defined by a linearized characteristic velocity. An alternative way to ensure conservation consists of transporting the entire control volume backward in time and subsequently integrating the advected quantity, as proposed in \cite{SL_Lipscomb_2005,bergami2023}. A conservative first-order finite volume scheme employing a semi-Lagrangian discretization of the advection terms was presented in \cite{SL_Lin_1996}, where diffusive effects are neglected. The accuracy and stability of Godunov-type solvers capable of operating with arbitrarily large time steps were analyzed in \cite{SL_Morton_1998} for general scalar conservation laws. In \cite{SL_Leonard_1996}, a family of advection schemes combining non-conservative and conservative advection formulations was developed to mitigate the errors introduced by operator-splitting approaches. Among the existing semi-Lagrangian methods, flux-based characteristic schemes \cite{SL_Frlkovic_2002} are also noteworthy, as they merge the conservation properties of finite volume discretizations with the large time-step capability inherent to semi-Lagrangian advection treatments.

Recently in \cite{SLIMEX}, a conservative $\SL$ method has been designed inspired by the approach proposed in \cite{QuiShu2011}. The new $\SL$ scheme was benchmarked for the nonlinear advection-diffusion equation and the one-dimensional shallow water model, combining an implicit-explicit time marching technique with a $\SL$ discretization. The method relies on integrating the governing equations over a space-time control volume that is built upon the backward trajectory of a pseudo-characteristic for each grid point. By applying Gauss theorem, the flux across the space-time surface drawn by the pseudo-characteristic is shown to vanish. In principle, this formulation achieves unconditional stability. However, for the benchmarks presented in \cite{QuiShu2011,SLIMEX}, the maximum $\CFL$ number employed was $\CFL=3$. For the results concerning the nonlinear Burgers equation shown in \cite{ConsSL2021}, the maximum attained $\CFL$ number was $\CFL=2.6$. In this work, instead, we employ space-time control volumes defined by the real backward characteristics (i.e. not the pseudo-characteristics) and deal with the consequent fact that, for nonlinear problems, the integral across the space-time surface associated with the characteristic might no longer vanish. In this way we obtain a novel conservative and unconditionally stable ($\CFL \gg 1$) semi-Lagrangian method for advection and diffusion problems. More specifically, computing the foot of the characteristic requires solving a nonlinear equation but enhances robustness and stability. For linear problems, no nonlinear equation has to be solved and, also in the novel method, the flux across the space-time surface of the characteristic automatically vanishes. Furthermore, for the first time, the diffusion terms are directly incorporated into a conservative semi-Lagrangian formulation at the aid of a novel characteristic-based Crank-Nicolson scheme. Our approach is therefore different from the discretization of parabolic terms in divergence form proposed in \cite{bonaventura2014siam} in the context of non-conservative semi-Lagrangian schemes, which makes use of a diffusive displacement that is directly accounted for in the particle trajectory. Likewise, stochastic streamlines are not used in our method, as instead done for linear advection-diffusion-reactions problems in \cite{bonaventura2021advdiff}. Eventually, we obtain a fully discrete space-time second order conservative finite difference scheme capable of handling nonlinear advection-diffusion problems with $\CFL$ numbers up to $100$.

The remainder of this article is organized as follows. Section \ref{sec:advection} introduces the advection equation and presents conservative semi-Lagrangian ($\SL$) schemes for both linear and nonlinear cases. In Section \ref{sec:advdiff}, diffusion terms are incorporated through a novel $\SL$ formulation, described in full detail. Section \ref{sec:numTests} validates the proposed $\SL$ scheme using a broad set of test problems covering linear and nonlinear advection as well as advection-diffusion scenarios. Finally, conclusions and directions for future research are provided in Section \ref{sec:concl}.

\section{Advection} \label{sec:advection}

Let us consider the one-dimensional advection equation of a scalar quantity $q=q(x,t)$ over a velocity field $u=u(x,t)$ in its conservative form:
\begin{equation}\label{eq:adv1d}
\begin{cases}
	\pder{t}{q} + \pder{x}{(uq)} = 0, & \text{for } x,t \in \R \times (0,T],\\
    q(x,0) = q_0(x), & \text{for } x \in \R,
\end{cases}
\end{equation}
where $t\in\R^+_0$ is the time, bounded in the interval $(0,T]$, and $x\in\R$ denotes the spatial coordinate. The model problem \eqref{eq:adv1d} is linear and, since the velocity $u$ is prescribed, its solution can be easily determined by solving the trajectory equation for the characteristics, that is
\begin{equation}\label{eq:odeCharLinear}
	\frac{\intd y}{\intd s} = -u(y,t-s), \quad s \in (0,t], \quad y(0) = x.
\end{equation}
More in general, we are interested in nonlinear hyperbolic conservation laws of the form
\begin{equation}\label{eq:conslaw1d}
\begin{cases}
	\pder{t}{q} + \pder{x}{f} = 0, & \text{for } x,t \in \R \times (0,T],\\
    q(x,0) = q_0(x), & \text{for } x \in \R,
\end{cases}
\end{equation}
where $f=f(q(x,t))$ is the flux function and for which the trajectory equation becomes
\begin{equation}\label{eq:odeCharF}
	\frac{\intd y}{\intd s} = -f'(q(y,t-s)), \quad s \in (0,t], \quad y(0) = x,
\end{equation}
with $f'=\partial f/\partial q$. As an example, for the inviscid Burgers' equation one gets $f=q^2/2$ and $f'=q$.

The complexity of the problem, and in particular its nonlinearity, is determined by the dependency of the right-hand-side in \eqref{eq:odeCharF} on the solution $q$. In such cases, even if the solution is initially smooth, the characteristics are modified by the solution itself, possibly determining the formation of shock waves. Capturing the correct speed of these waves is of paramount importance when dealing with conservation laws, and thus for the numerical scheme used to approximate their solution. Therefore, in what follows, we aim at constructing a conservative $\SL$ scheme ($\SCOUT$) for the general problem \eqref{eq:conslaw1d}.\\

Referring to \eqref{eq:conslaw1d}, we first discretize the problem in time. Let $\dt>0$ be a time step, $\tn = n\dt$ a uniform time grid with $n = 0, \ldots, N_T$, where $N_T = \lceil \frac{T}{\dt} \rceil$. We start by considering the semi-discrete form
\begin{equation}\label{eq:semidiscr1d}
    q(x,\tnn) = q(x,\tn) - \PDER{x}{}\int_{\tn}^{\tnn} f(q(x,t)) \, \intd t,
\end{equation}
where we keep the time integral notation for the flux term, as it will be useful in the sequel. From \eqref{eq:semidiscr1d}, it follows that updating the solution $q$ lies in the computation, and subsequent derivation, of a flux function $H$, which it is defined as
\begin{equation}\label{eq:Hflux}
    H(x) \coloneqq \int_{\tn}^{\tnn} f(x,t,q(x,t)) \, \intd t.
\end{equation}

Let us now consider a computational domain $X=[a,b]$ and let us introduce a space grid of uniform size $\dx$. Given $N_x$ the total number of cells used to discretize $X$, each cell has a constant spacing $\dx=(b-a)/N_x$. For a cell, its barycenter is denoted with $x_i$ and it is centered within the interval $[x_{i-1/2}, x_{i+1/2}]$ with the interfaces given by $x_{i\pm1/2}=x_i\pm \dx/2$. We adopt a cell-centred space discretization, so that the conserved quantity $q$ is defined at the cell barycentres $x_i$, namely $q^n_i=q(x_i,\tn)$. With this setting, a conservative finite difference scheme can be derived from \eqref{eq:semidiscr1d} and \eqref{eq:Hflux} as
\begin{equation}\label{eq:numSchemeGeneric}
    q^{n+1}_i = q^n_i - (D_x \hat{H})_i,
\end{equation}
where $D_x$ is a discrete version of the differential operator $\pder{x}{}$ and $\hat{H}$ is a suitable numerical approximation for the function $H$. In particular, resorting to a second order central-difference approximation for $D_x$, the scheme takes the general conservative form
\begin{equation}\label{eq:numScheme}
    q^{n+1}_i = q^n_i - \frac{1}{\dx} \left( \hat{H}_{\ipm} - \hat{H}_{\imm}\right),
\end{equation}
The remaining question is how to obtain the numerical fluxes $\{ \hat{H}_{\ipm}\}_{i=0}^{N_x}$. For the sake of clarity, we will detail this part distinguishing between the linear case and the general nonlinear one.

\subsection{Linear advection}\label{ssec:linearCase}

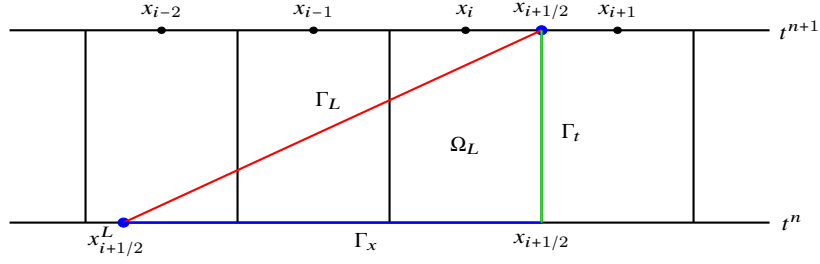
\begin{figure}[!h]
	\centering
	\begin{tikzpicture}[xscale=2.,yscale=1.7]
    \tikzmath{\x=1; \y=1.5; \xL=-1.75; \yL=0; \a=-2.5; \b=2.5;} 

    \draw[thick] (\a,\y) -- (\b,\y);

    \draw[thick] (\a,\yL) -- (\b,\yL);

    \draw[thick] (0,\y) -- (0,\yL);
    \draw[thick] (+1,\y) -- (+1,\yL);
    \draw[thick] (+2,\y) -- (+2,\yL);
    \draw[thick] (-1,\y) -- (-1,\yL);
    \draw[thick] (-2,\y) -- (-2,\yL);

    \filldraw[blue] (\x,\y) circle (1pt);
    \filldraw[blue] (\xL,\yL) circle (1pt);
    \node at (\x,\y+0.15) {$\xipm$};
    \node at (\xL-0.05,\yL-0.15) {$x^L_{\ipm}$};
    \node at (\x,\yL-0.15) {$\xipm$};

    \node at (\x-0.5,\y+0.15) {$x_i$};
    \node at (\x-1.5,\y+0.15) {$x_{i-1}$};
    \node at (\x-2.5,\y+0.15) {$x_{i-2}$};
    \node at (\x+0.5,\y+0.15) {$x_{i+1}$};

    \node at (\x-0.5,\y-0.9) {$\domainlag$};
    \node at (\x-1.4,\y-0.55) {$\Gamma_L$};
    \node at (\x-1.15,\yL-0.15) {$\Gamma_x$};
    \node at (\x+0.2,\y-0.8) {$\Gamma_t$};

    \filldraw[black] (\x-0.5,\y) circle (0.7pt);
    \filldraw[black] (\x-1.5,\y) circle (0.7pt);
    \filldraw[black] (\x-2.5,\y) circle (0.7pt);
    \filldraw[black] (\x+0.5,\y) circle (0.7pt);

    \draw[red,thick] (\x,\y) to (\xL,\yL);
    \draw[blue,thick] (\xL,\yL) to (\x,\yL);
    \draw[green,thick] (\x,\yL) to (\x,\y);

    \node at (\b+0.2,\y) {$\tnn$};
    \node at (\b+0.15,\yL) {$\tn$};

\end{tikzpicture}
	\caption{Space-time domain $\domainlag$ used to derive the conservative $\SL$ scheme in the linear case with constant velocity.}
	\label{fig:OmegaL}
\end{figure}

Let us first consider the simplified problem \eqref{eq:adv1d} with constant velocity $u(x,t)=\Bar{u}$, for which the trajectory equation takes the form
\begin{equation}\label{eq:odeCharConst}
    \frac{\intd y}{\intd s} = -\Bar{u},
\end{equation}
with $s=T-t$, which states that the solution of \eqref{eq:adv1d} is advected along parallel straight lines having slope $\Bar{u}$. By solving \eqref{eq:odeCharConst}, we can associate to each point $\xipm$ the Lagrangian trajectory $\Gamma_L=\{(s,y(s):s\in[0,\dt]\}$, its foot of the characteristic $x^L_{\ipm}=y(\dt)$, and a corresponding space-time domain $\domainlag$ that is the region bounded by the segments $\Gamma_x=[x^L_{\ipm},\xipm]$, $\Gamma_t=[\tn,\tnn]$, and the Lagrangian trajectory $\Gamma_L$, as shown in Fig.~\ref{fig:OmegaL}.

To enforce conservation, we consider the integral of the PDE \eqref{eq:adv1d} over the domain $\domainlag$ and define the space-time divergence operator $\nabla_{x,t} \coloneqq (\pder{x}{},\pder{t}{})$, so that the use of Gauss theorem in the space-time framework gives the equality
\begin{equation}\label{eq:PDEintegralConst}
\begin{aligned}
    0 = \int_{\domainlag}\pder{t}{q}+\pder{x}{\Bar{u}q}\,\intd x\, \intd t &= \int_{\domainlag} \nabla_{x,t} \cdot (\Bar{u}q,q) \,\intd x\, \intd t\\
    &= \int_{\Gamma_L} (\Bar{u}q,q) \cdot \vec{n}_L \, \intd s + \int_{\Gamma_x} (\Bar{u}q,q) \cdot \vec{n}_x \, \intd x + \int_{\Gamma_t} (\Bar{u}q,q) \cdot \vec{n}_t \, \intd t.
\end{aligned}
\end{equation}
By construction, in the $(x,t)$ reference system, since the tangent vector to $\Gamma_L$ is given by the velocity of the trajectory, i.e. $\vec{t}_L=(-\Bar{u},-1)$, the corresponding normal vector results to be
\begin{equation}\label{eq:normalLConst}
\vec{n}_L=\lambda\, \beta\, (1,-\Bar{u}),
\end{equation}
where $\lambda=\sign(\Bar{u})$ and $\beta = 1/\sqrt{1+\Bar{u}^2}$. Therefore the integral along $\Gamma_L$ in \eqref{eq:PDEintegralConst} vanishes. According to the sign, we also get $\vec{n}_x=(0,-1)$ and $\vec{n}_t=\lambda\,(1,0)$, which, rearranging the terms in \eqref{eq:PDEintegralConst}, leads to
\begin{equation}\label{eq:fluxItegralExaConst}
    \int_{\tn}^{\tnn} \Bar{u}q
    \, \intd t = \lambda\int_{x^L_{\ipm}}^{\xipm} q\, \intd x.
\end{equation}
As a consequence, from \eqref{eq:fluxItegralExaConst}, it follows that the flux function ${H}(x)$ can be computed just relying on spatial information at time $\tn$ as
\begin{equation}\label{eq:HcomputationConst}
    {H}(\xipm)=\lambda\int_{x^L_{\ipm}}^{\xipm} q(x,\tn) \, \intd x.
\end{equation}
Note that, when applying Gauss theorem in \eqref{eq:PDEintegralConst}, our convention is to trace the Lagrangian path backward, and consequently proceed along the other edges, $\Gamma_x$ and $\Gamma_t$. The same holds true in the subsequent cases, regardless of the position of the foot $x^L_{\ipm}$ with respect to the node $\xipm$. 

The computation of the numerical fluxes $\hat{H}_{\ipm}$ and $\hat{H}_{\imm}$, appearing in \eqref{eq:numScheme}, requires a proper numerical approximation of the integral \eqref{eq:HcomputationConst} that will be detailed later in Sec.~\ref{ssec:numApprox}.

\paragraph{Space- and time-dependent velocity}
The discretization \eqref{eq:HcomputationConst} for the flux function $H$ follows analogously when the advection velocity $u(x,t)$ is not constant, but explicitly depends on space and time. Assuming that the sign of $u$ remains constant in the interval $[\tn,\tnn]$, the only difference concerns the computation of the point $x^L_{\ipm}$. Since the characteristics are no longer straight lines, solving the trajectory equation might not be straightforward and numerical integration should be properly employed, as discussed later in Sec.~\ref{ssec:numApprox}.

\subsection{Nonlinear case}\label{ssec:nonlinearCase}
Let us now consider the more complex case in which, referring to \eqref{eq:conslaw1d}, the flux is given by $f=f(q(x,t))$ and the solution is advected according to \eqref{eq:odeCharF} by $f'=f'(q(x,t))$. In this case, the characteristics are still straight lines, but have a slope which depends on the solution itself. Similarly to the previous case, and referring again to Fig.~\ref{fig:OmegaL}, we can provide a similar construction of the domain $\domainlag$ and get the equality
\begin{equation}\label{eq:PDEintegral}
    \int_{\Gamma_L} (f,q) \cdot \vec{n}_L \, \intd s + \int_{\Gamma_x} (f,q) \cdot \vec{n}_x \, \intd x + \int_{\Gamma_t} (f,q) \cdot \vec{n}_t \, \intd t = 0,
\end{equation}
from PDE system \eqref{eq:conslaw1d}. For the sake of clarity, we now assume that $f'(q(\xipm,t))$ does not change its sign in the interval $[\tn,\tnn]$, so that we can define
\begin{equation}
    \lambda\coloneqq\sign(f'(q(\xipm,\tnn))).
\end{equation}
To lighten notation, in what follows we omit the dependence of the function $f'$ on $s=T-t$, as well as that of the other functions on the variables of the problem, since the value of the solution is constant along the characteristics.
With this setting, we immediately get $\vec{n}_x = (0,-1)$ and $\vec{n}_t = \lambda\,(1,0)$, while tracing backward the characteristic we have the tangent vector
\begin{equation}\label{eq:normalL}
    \vec{t}_L = (-f',-1) \quad \text{and thus} \quad \vec{n}_L = \lambda\,\beta\,(-1,f'),
\end{equation}
where $\beta=1/\sqrt{1+(f')^2}$.

The crucial aspect now lies in observing that, in general, compared to the linear case, the integral on $\Gamma_L$ does not vanish any more. In fact, from \eqref{eq:PDEintegral} and \eqref{eq:normalL} we get the equality
\begin{equation}\label{eq:fluxItegralExa}
    \int_{\tn}^{\tnn} f \,\intd t = \lambda\int_{x^L_{\ipm}}^{\xipm} q\, \intd x + \beta \int_{(\xipm,\tnn)}^{(x^L_{\ipm},\tn)} (f-qf')\, \intd s,
\end{equation} 
which differs from \eqref{eq:fluxItegralExaConst} due to the last term. As an example, for the nonlinear Burgers' equation one gets $f-qf'=-q^2/2$, while for the linear case one recovers $f=uq$, $f'=u$ and hence $f-qf'=uq-uq=0$. In \cite{QuiShu2011}, instead of the characteristic speed $f'$, the trajectory equation is governed by a pseudo-characteristic speed inferred from the conservative formulation as $f(q)/q$, so the last integral in \eqref{eq:fluxItegralExa} vanishes. In the Burgers' example, the pseudo-characteristic speed employed in \cite{QuiShu2011} is $q/2$, while the characteristic speed is $f'=q$.

The computation of $\hat{H}_{\ipm}$ in the nonlinear case is thus based on the evaluation of two integral terms: the first one only needs spatial information at time $\tn$ whereas the last one requires an integral computation along the Lagrangian trajectory. However, once the foot $x^L_{\ipm}$ is computed, since the solution is constant along the characteristic, one can write
\begin{equation}\label{eq:HcomputationNonlinear}
    H(x_{\ipm}) = \lambda\int_{x^L_{\ipm}}^{\xipm} q \,\intd x + \dt \Big[(f-qf')\big|_{(x^L_{\ipm},\tn)}\Big],
\end{equation}
where the arclength of the curve $\Gamma_L$ simplifies with $\beta$, thus reducing the last term to a multiplication between the time step $\dt$ and the integrand, evaluated at the foot $(x^L_{\ipm},\tn)$.

At this point, it is worth noting that no numerical approximation has been involved yet. In what follows, the fully discrete numerical scheme will be derived by introducing suitable techniques.

\subsection{Numerical approximation of the flux function}\label{ssec:numApprox}
The discretization $\hat{H}_{\ipm}$ of the flux functions $H_{\ipm}$ given by \eqref{eq:HcomputationConst} or \eqref{eq:HcomputationNonlinear}, to be used in the conservative scheme \eqref{eq:numScheme}, requires three main ingredients:
\begin{itemize}
    \item the computation of the feet $\{ x^L_{\ipm}\}_{i=0}^{N_X}$, by integrating the trajectory equation \eqref{eq:odeCharF};
    \item a quadrature rule for the approximation of the space integral in \eqref{eq:HcomputationConst} or \eqref{eq:HcomputationNonlinear}, which for convenience we define as
    \begin{equation}\label{eq:spatialIntegral}
        I_{\ipm} \coloneqq \int_{x^L_{\ipm}}^{\xipm} q(x,\tn)\, \intd x;
    \end{equation}
    \item a spatial reconstruction operator $\rec[Q^n]$ for the solution $q$ at time $\tn$, where $Q^n$ denotes the set of point values $\{ q_i\}_{i=1}^{N_x}$ at time $\tn$.
\end{itemize}
We point out that, in this work, we aim at designing a second order accurate scheme in both space and time.

\paragraph{Computation of the Lagrangian feet.} In the linear case, integrating \eqref{eq:odeCharLinear} is straightforward since the velocity $u(x,t)$ is prescribed. In particular, when $u(x,t)=\Bar{u}$ is constant, we can exactly compute
\begin{equation}\label{eq:RKconst}
    x^L_{\ipm} = \xipm - \Bar{u}\dt.
\end{equation}
Otherwise, when the velocity explicitly depends on space and time, a suitable $m$-stage Runge Kutta (RK) method can be applied as follows
\begin{equation}\label{eq:RK}
\begin{aligned}
 x^L_{\ipm} &= \xipm + \dt\sum_{r=1}^{m} b_r K_r,\\
 K_r &= -u\,(Y_r,t_n+(1-c_r)\dt),\\
 Y_r &= x_{\ipm}+\dt\sum_{j=1}^{r-1}A_{rj}K_j,
\end{aligned}
\end{equation}
where $b_r$,$c_r$,$A_{rj}$ are the coefficients of the Butcher tableau defining the RK method used to solve \eqref{eq:odeCharLinear}. In particular, we resort to Heun's method
\begin{equation}\label{eq:heun}
\begin{array}
{c|cc}
0\\
1 & 1\\
\hline
& \nicefrac{1}{2} &\nicefrac{1}{2}
\end{array}
\end{equation}
to get second order accuracy. Note that, both in \eqref{eq:RKconst} and \eqref{eq:RK}, the signs are consistent with backward integration.

On the other hand, for  nonlinear problems the velocity depends on the solution itself and it is not given \textit{a priori}. Let us first integrate \eqref{eq:odeCharF} in the time interval $[0,\dt]$ in order to get the equality
\begin{equation}\label{eq:odeIntegration}
    x^L_{\ipm} = \xipm - \int_{0}^{\dt} f'(q(y,\tnn-s))\, \intd s,
\end{equation}
to which a suitable quadrature rule can be applied. For our purposes, employing the trapezoidal rule leads to
\begin{equation}\label{eq:odeIntegrationTrapezoidal}
    x^L_{\ipm} = \xipm - \frac{1}{2}\left[ f'(q(x^L_{\ipm},\tn)) + f'(q(\xipm,\tnn)) \right]\dt,
\end{equation}
which simplifies to
\begin{equation}\label{eq:odeIntegrationTrapezoidalSimplified}
    x^L_{\ipm} = \xipm - f'(q(x^L_{\ipm},\tn))\dt,
\end{equation}
by using the fact that $q$ is constant along the characteristics, thus enforcing $q(\xipm,\tnn)=q(x^L_{\ipm},\tn)$. In order to compute the solution $x^L_{\ipm}$ of the nonlinear equation \eqref{eq:odeIntegrationTrapezoidalSimplified}, we adopt Newton's method, setting as initial guess
\begin{equation}\label{eq:NewtonInitialGuess}
    x_0 = \xipm - f'(\rec[Q^n](\xipm))\dt.
\end{equation}
In practice, in order to guarantee the robustness of the numerical scheme, we also employ a parachute bisection method when Newton's algorithm fails to converge.

\paragraph{Quadrature rule.} The numerical approximation of the spatial integral $I_{\ipm}$ given by \eqref{eq:spatialIntegral}, to be used in \eqref{eq:HcomputationConst} or \eqref{eq:HcomputationNonlinear}, relies on a piecewise polynomial representation of the solution $Q^n$ and on a Gauss-Legendre quadrature rule with $N_K$ nodes. Let us remark that the path defined by $\Gamma_x$ may cross several cells. 

We denote by $j^*$ the index of the cell in which the characteristic foot $x^L_{\ipm}$ lies, and by $\dx^*$ the length $\lvert x^L_{\ipm} - x_{j^*+1/2} \rvert$. Then, we can approximate the integral $I_{\ipm}$ as follows.
\begin{itemize}
    \item If $x^L_{\ipm}\leq\xipm$ then
    \begin{equation}\label{eq:spatialIntegralApproxLeft}
    \begin{aligned}
    I_{\ipm} &= \int_{x^L_{\ipm}}^{x_{j^*+1/2}} q\, \intd x + \sum\limits_{j=j^*+1}^i \int_{x_{j-1/2}}^{x_{j+1/2}} q\, \intd x \\
    &\approx \dx^*\sum\limits_{k=1}^{N_K} w_k q(x^L_{\ipm}+\xi_k\dx^*,\tn) + \sum\limits_{j=j^*+1}^i \dx\sum\limits_{k=1}^{N_K} w_k q(x_{j-1/2}+\xi_k\dx,\tn)\\
    &\approx \dx^*\sum\limits_{k=1}^{N_K} w_k \rec[Q^n](x^L_{\ipm}+\xi_k\dx^*) + \sum\limits_{j=j^*+1}^i \dx\sum\limits_{k=1}^{N_K} w_k \rec[Q^n](x_{j-1/2}+\xi_k\dx).
    \end{aligned}
    \end{equation}
    \item If $x^L_{\ipm}>\xipm$ then
    \begin{equation}\label{eq:spatialIntegralApproxRight}
    \begin{aligned}
    I_{\ipm} &= \int_{x^L_{\ipm}}^{x_{j^*-1/2}} q\, \intd x + \sum\limits_{j=j^*-1}^{i+1} \int_{x_{j+1/2}}^{x_{j-1/2}} q\, \intd x \\
    &\approx (\dx-\dx^*)\sum\limits_{k=1}^{N_K} w_k q(x^L_{\ipm}-\xi_k\dx^*,\tn) + \sum\limits_{j=j^*-1}^{i+1} \dx\sum\limits_{k=1}^{N_K} w_k q(x_{j+1/2}-\xi_k\dx,\tn)\\
    &\approx (\dx-\dx^*)\sum\limits_{k=1}^{N_K} w_k \rec[Q^n](x^L_{\ipm}-\xi_k\dx^*) + \sum\limits_{j=j^*-1}^{i+1} \dx\sum\limits_{k=1}^{N_K} w_k \rec[Q^n](x_{j+1/2}-\xi_k\dx).
    \end{aligned}
    \end{equation}
\end{itemize}
Here, $\xi_k$ and $w_k$ are the nodes and the weights of the quadrature rule, respectively. We point out that, in our work, a simple midpoint rule is employed in order to achieve second order of accuracy, thus $N_K=1$ with $\xi_1=0.5$ and $w_1=1$.

\paragraph{Spatial reconstruction} A piecewise linear reconstruction is implemented for the approximation of the spatial values $Q^n$, which, for a generic cell $i$, is given by
\begin{equation}\label{eq:linearReconstruction}
    \rec[Q^n]_i(x) = q_i^n + \sigma_i^n (x-x_i)/\dx.
\end{equation}
The slope is simply given by
\begin{equation}\label{eq:linearReconstructionSlopes}
    \sigma_i^n = \frac{1}{2}\left(q_{i+1}-q_{i-1}\right),
\end{equation}
which could be bounded by a TVD limiter, e.g. minmod. In all the test cases shown in Sec.~\ref{sec:numTests} we simply use the central reconstruction \eqref{eq:linearReconstructionSlopes} with no limiter.

\section{Advection-diffusion} \label{sec:advdiff}
We now turn our attention to diffusion processes. In order to simplify the computation and without loss of generality we consider here a linear advective contribution. The extension to nonlinear advection can easily be obtained following the details reported in Sec.~\ref{ssec:nonlinearCase}. Starting from \eqref{eq:adv1d}, the mathematical model now writes
\begin{equation}\label{eq:advdiff1d}
	\begin{cases}
		\pder{t}{q} + \pder{x}{(uq)} = \pder{x}{\nu\pder{x}{q}}, & \text{for } x,t \in \R \times (0,T],\\
		q(x,0) = q_0(x), & \text{for } x \in \R,
	\end{cases}
\end{equation}
where $\nu\geq 0 $ is a constant diffusion coefficient. 
Integration of \eqref{eq:advdiff1d} over the time interval $[\tn,\tnn]$ leads to
\begin{eqnarray} \label{eq:advdiff_divint}
	\int_{\tn}^{\tnn} \pder{t}{q} \,\intd t &=&  -\int_{\tn}^{\tnn}\pder{x}{\left(uq-\nu\pder{x}{q}\right)}\,\intd t, \nonumber \\
    q(x,\tnn)-q(x,t^{n})&=&- \pder{x}{\left(\int_{\tn}^{\tnn}{uq \,\intd t}-\int_{\tn}^{\tnn}{\nu \pder{x}{q} \,\intd t}\right)},
\end{eqnarray}
which does not involve any numerical approximation yet. In order to evaluate the time integrals in \eqref{eq:advdiff_divint} with a conservative semi-Lagrangian approach, the governing PDE \eqref{eq:advdiff_divint} is integrated over the domain $\domainlag$ originating from a generic point $x$, hence yielding
\textcolor{black}{
\begin{eqnarray}
\label{eq:advdiff_stint}
	    \int_{\domainlag}\pder{t}{q}+\pder{x}{uq}\,\intd x\, \intd t &=& \int_{\domainlag} \pder{x}{(\nu \pder{x}{q})} \,\intd x \,\intd t, \nonumber \\
	    \int_{\tn}^{\tnn}uq \,\intd t -\lambda\int_{x^L}^{x}q \,\intd \xi
	      &=& \int_{\tn}^{\tnn} \int_{y^L(s)}^{x}\pder{x}{(\nu \pder{x}{q})} \,\intd x \,\intd t, 
\end{eqnarray}
where $y^L(t)$ indicates the characteristic location starting from $x$ at a generic time $t \in [t^n,t^{n+1}]$, hence $y^L(\tnn)=x$ and $y^L(\tn)=x^L$. Let us remark that in \eqref{eq:advdiff_stint} we keep the notation $\lambda=\sign(u)$, assuming that it remains constant in the interval $[\tn,\tnn]$.}
Let us now introduce the abbreviation 
\begin{equation}
	\label{eqn.alpha}
	\alpha(x,t):=\pder{x}{q(x,t)}
\end{equation} 
for the space derivative of $q$. Application of the fundamental theorem of calculus over the last integral in \eqref{eq:advdiff_stint} leads to
\begin{eqnarray} \label{eq:advdiff_stint2}
	\int_{\tn}^{\tnn}uq \,\intd t&=&\lambda \int_{x^L}^{x}q(\xi,\tn) \,\intd \xi  +\nu \int_{\tn}^{\tnn}  \alpha(x,t)-\alpha(y^L(t),t) \,\intd t
\end{eqnarray}
Substituting \eqref{eq:advdiff_stint2} into \eqref{eq:advdiff_divint} we get
\begin{eqnarray}
q(x,\tnn)-q(x,\tn) &=&-\partial_x\left(\int_{\tn}^{\tnn} uq\, \intd t-\int_{\tn}^{\tnn}\nu \alpha\,\intd t \right) \nonumber \\
	&=&-\partial_x\left(\lambda\,\,\int_{x^L}^{x}q(\xi,\tn) \,\intd \xi-\nu \int_{\tn}^{\tnn}\alpha(y^L(t),t)\, \intd t\right),
\end{eqnarray}
that is still exact in the current form and, with the usual notation for the function $H(x)$ according to \eqref{eq:Hflux}, yields the semi-discrete scheme
\begin{equation}
	\label{eqn.advdiff_sd}
    q(x,\tnn) = q(x,\tn) - \partial_x \left( H(x) - \nu \int_{\tn}^{\tnn}\alpha(y^L(t),t)\, \intd t\right).
\end{equation}

The last time integral in \eqref{eqn.advdiff_sd} is approximated by means of the trapezoidal rule to achieve a second order of accuracy, thus one gets
\textcolor{black}{
\begin{eqnarray}
	\int_{\tn}^{\tnn}\alpha(y^L(t),t)\,\intd t &=& \frac{\dt}{2}\left[\alpha(y^L(\tnn),\tnn)+\alpha(y^L(\tn),\tn)\right] + \Ogrande(\dt^2),
\end{eqnarray}
that separates the fully implicit contribution $\alpha(y^L(\tnn),\tnn)=\alpha(x,\tnn)$ and the explicit one, which is computed at the foot of the characteristic.}

\subsection{Semi-implicit scheme} \label{ssec:siadvdif}
The fully discrete semi-implicit scheme is obtained relying on a finite difference approximation in space. Therefore, the spatial derivative in \eqref{eqn.advdiff_sd} is approximated by the central finite difference operator already used in \eqref{eq:numScheme} where the numerical flux functions of $\hat{H}_{i\pm1/2}$ are the same as detailed in Sec.~\ref{ssec:numApprox}. The scheme then writes
\begin{eqnarray}\label{eq:discadvdif}
	q_i^{n+1} = q_i^n-\frac{1}{\Delta x}\left(\hat{H}_{\ipm}-\hat{H}_{\imm} \right)+\frac{\nu \Delta t}{2\Delta x}\left[\hat\alpha^{n+1}_{\ipm}-\hat\alpha^{n+1}_{\imm}\right]+\frac{\nu \Delta t}{2 \Delta x}\left[\alpha^n(x^L_{\ipm})-\alpha^n(x^L_{\imm}) \right], \nonumber \\
\end{eqnarray}
where the approximations $\hat\alpha^{n+1}_{i\pm1/2}$ are given by
\begin{equation}\label{eq:alphaapprox}
\alpha(\xipm,\tnn)\approx\hat\alpha^{n+1}_{\ipm}=\frac{q_{i+1}^{n+1}-q_i^{n+1}}{\Delta x}, \qquad \alpha(\ximm,\tnn)\approx\hat\alpha^{n+1}_{\imm}=\frac{q_{i}^{n+1}-q_{i-1}^{n+1}}{\Delta x}.
\end{equation}
With these definitions, the fully discrete scheme can be compactly rewritten as
\begin{eqnarray}\label{eq:linsys}
	-\frac{\nu \Delta t}{2\Delta x^2}q_{i-1}^{n+1}+\left(1+\frac{\nu \Delta t}{\Delta x^2}\right)q_i^{n+1}-\frac{\nu \Delta t}{2\Delta x^2}q_{i+1}^{n+1} = b_i^n,
\end{eqnarray}
with the explicit right-hand-side contribution
\begin{eqnarray}\label{eq:linsysb}
	b_i^n=q_i^n-\frac{1}{\Delta x}\left(\hat{H}_{\ipm}-\hat{H}_{\imm} \right)+
	\frac{\nu \Delta t}{2 \Delta x}\left[\alpha^n(x^L_{\ipm})-\alpha^n(x^L_{\imm}) \right].
\end{eqnarray}
The tri-diagonal linear system \eqref{eq:linsys} is efficiently solved using the Thomas algorithm to obtain the new solution $q_i^{n+1}$. If no advection is present, e.g. when $u=0$, or the advection effects are negligible compared to the viscous contribution, the scheme \eqref{eq:discadvdif} reduces to the classical second order Crank-Nicolson approximation of the diffusion terms. However, when advection processes are important, the time integral of the diffusion terms is computed along the characteristics, that is fundamental in order to achieve second order of accuracy. This is where the novel semi-Lagrangian Crank-Nicolson discretization plays a key role. Indeed, the diffusion fluxes are evaluated as integrals along the Lagrangian trajectories as investigated in \cite{SLIMEX}. In other words, when the time step is large ($\text{CFL} \gg 1$), or the diffusion process is not negligible with respect to the advective one, the two phenomena cannot be decoupled and a consistent discretization of the parabolic terms is essential to achieve the proper accuracy. We refer to this scheme with semi-Lagrangian Crank-Nicolson. We point out that this scheme shares with the ones of \cite{bonaventura2014siam,bonaventura2021advdiff} the idea of approximating the diffusion term at the foot of the characteristic; there, however, this term is reinterpreted as Brownian motion and approximated by evaluating the reconstruction at time $t^n$ at a distance $\pm\sqrt{\dt}$ from the foot, obtaining an explicit but not conservative scheme. Here, we propose a characteristic-based Crank-Nicolson discretization that directly accounts for the diffusion term along the characteristics in a fully conservative formulation for nonlinear problems.

\section{Numerical results} \label{sec:numTests}

In this section, numerical tests for one-dimensional linear and nonlinear hyperbolic equations are presented in order to validate the new numerical method. In particular, we focus on the expected accuracy of the $\SCOUT$ scheme under very large Courant numbers, considering the $\L^1$-norm of the error and the conservation of mass, by computing the quantity
\begin{equation}
    \Delta q = \Big\lvert \left\lVert Q^{N_T} \right\rVert_1 - \left\lVert Q^0 \right\rVert_1 \Big\rvert.    
\end{equation}
To this aim, the novel conservative scheme is compared with a classical $\SL$ one based on the same Runge-Kutta time integrator and on the same spatial reconstruction. In this comparison, we expect the classical $\SL$ scheme to achieve first order of accuracy, according to the analysis carried out in \cite{Carlini2013}, and, in addition, to be non-conservative. In all figures we refer to $\SCOUT$ and $\SL$ to precisely address these two schemes. Furthermore, for benchmarks involving advection-diffusion terms, we refer to the novel semi-Lagrangian Crank-Nicolson scheme proposed in Sec.~\ref{sec:advdiff} and to the classical Crank-Nicolson method with $\SCOUTSLCN$ and $\SCOUTCN$, respectively. 

We point out that many of the tests reported here are inspired by the ones presented in \cite{ConsSL2021}, as they are well apt at verifying the effectiveness, robustness and accuracy of the $\SCOUT$ method. 

We recall that $N_x$ denotes the number of cells and $T$ is the final time. The time step is computed according to a classical CFL-type condition, so that
\begin{equation}\label{eq:CFLtimestep}
    \dt = \CFL \frac{\dx}{\max\limits_x \lvert f'(x,t,q(x,t))\rvert},
\end{equation}
where we can adopt very large values for the $\CFL$ number since the scheme is unconditionally stable. Some example of truly large $\CFL$ numbers will be reported in the tests to empirically show this important property exhibited by $\SCOUT$ scheme.

\paragraph{Test 1: Linear constant advection of Gaussian data.}

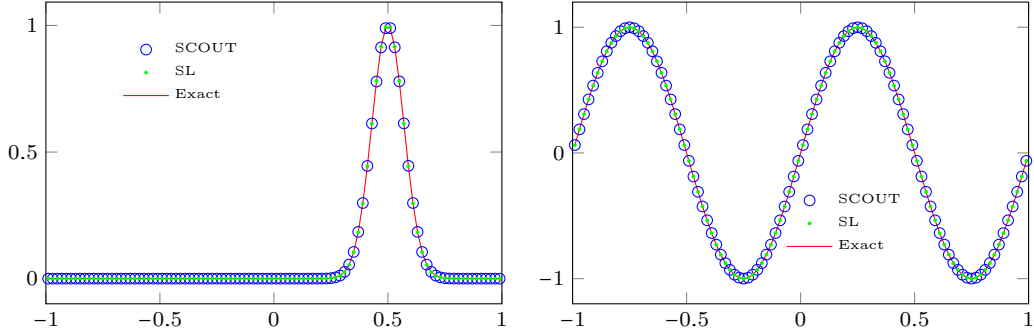
\begin{figure}[!h]
	\centering
	\begin{tikzpicture}

\pgfplotstableread{test1_convStep2.txt}{\graphData}

\begin{axis}[%
  width=6cm,height=4cm,scale only axis,
  xmin=-1,xmax=1,
  legend columns=1, 
  legend cell align={left},
  legend style={
    fill=none,
  	at={(0.45,0.9)},
  	draw=none,
    font=\tiny,
  }
]

\addplot[blue,mark=o,only marks] table[x=xcoor,y=solSLc] {\graphData};
\addlegendentry{$\SCOUT$};

\addplot[color=green,mark=*,only marks,mark size=0.5pt] table[x=xcoor,y=solSL] {\graphData};
\addlegendentry{$\SL$};

\addplot[color=red,mark=none] table[x=xcoor,y=solExa] {\graphData};
\addlegendentry{Exact};

\end{axis}
\end{tikzpicture}
	\begin{tikzpicture}

\pgfplotstableread{test3_convStep2.txt}{\graphData}

\begin{axis}[%
  width=6cm,height=4cm,scale only axis,
  xmin=-1,xmax=1,
  legend columns=1, 
  legend cell align={left},
  legend style={
    fill=none,
  	at={(0.75,0.4)},
  	draw=none,
    font=\tiny,
  }
]

\addplot[blue,mark=o,only marks] table[x=xcoor,y=solSLc] {\graphData};
\addlegendentry{$\SCOUT$};

\addplot[color=green,mark=*,only marks,mark size=0.5pt] table[x=xcoor,y=solSL] {\graphData};
\addlegendentry{$\SL$};

\addplot[color=red,mark=none] table[x=xcoor,y=solExa] {\graphData};
\addlegendentry{Exact};

\end{axis}
\end{tikzpicture}
	\caption{Constant linear advection case: comparison between the solution at the final time provided by the classical $\SL$ scheme, the novel conservative $\SCOUT$ scheme, and the exact one. Left: Test 1. Right: Test 2. In both tests, we set $N_x=100$, $T=1$, $\nu=0$ and $\CFL=10$.}
	\label{fig:graficoTest1}
\end{figure}

We first consider a benchmark test concerning the linear advection of a Gaussian initial profile
\begin{equation}
    q_0(x)=e^{-\left(\frac{x+0.5}{0.1}\right)^2},
\end{equation}
with constant velocity $\Bar{u}=1$. Then we are solving
\begin{equation}\label{eq:test1}
    \pder{t}{q} + \pder{x}{q} = 0,
\end{equation}
with $X=[-1,1]$, $T=1$, $\CFL=10$ and Dirichlet boundary conditions. The exact solution can be easily computed as $q(x,t)=q_0(x-t)$. In this example, since the trajectory equation \eqref{eq:odeCharConst} can be solved exactly, no temporal error exists and we expect to achieve second order of accuracy and to guarantee mass conservation with both schemes, namely $\SCOUT$ and $\SL$. In Tab.~\ref{tab:errorTest1} the $\L^1$-norm of the error and the corresponding order are reported, while the final solutions are depicted in the left panel of Fig.~\ref{fig:graficoTest1}. The mass error results in the order of the machine precision, as expected.

\begin{table}[!h]
	\centering
	\begin{tabular}{|c|c|c|c|c|c|c|}
		\hline \multirow{2}{*}{$N_x$} & \multicolumn{3}{c|}{$\SL$} & \multicolumn{3}{c|}{$\SCOUT$} \\
		\cline{2-7} & $\Delta q$ & $\L^1$-err & $L^{1}$-ord & $\Delta q$ & $\L^{1}$-err & $\L^1$-ord\\ \hline \hline
		$50$ & $0.00e+00$ & $2.74e-03$ & - & $0.00e+00$ & $2.74e-03$ & - \\ \hline
		$100$ & $1.02e-13$ & $1.11e-03$ & $1.31$ & $1.19e-13$ & $5.91e-04$ & $2.22$ \\ \hline
		$200$ & $3.50e-15$ & $1.93e-04$ & $2.52$ & $5.83e-15$ & $1.44e-04$ & $2.04$ \\ \hline
		$400$ & $2.50e-16$ & $4.01e-05$ & $2.27$ & $2.50e-16$ & $3.58e-05$ & $2.01$ \\ \hline
		$800$ & $5.55e-17$ & $9.27e-06$ & $2.11$ & $5.55e-17$ & $8.94e-06$ & $2.00$ \\ \hline
	\end{tabular}
	\caption{Test 1: mass conservation, errors and accuracy. Comparison between $\SCOUT$ and $\SL$ schemes.}
	\label{tab:errorTest1}
\end{table}

\paragraph{Test 2: Linear constant advection of periodic data.} The second benchmark test concerns the linear advection of a sinusoidal data
\begin{equation}
    q_0(x) = \sin(2\pi x)
\end{equation}
with constant velocity $\Bar{u}=-1$, in the domain $X=[-1,1]$, setting $T=1$, $\CFL=10$ and periodic boundary conditions. Similarly to the previous test, we can compute the exact solution as $q(x,t)=q_0(x+t)$ and verify that second order of accuracy and mass conservation are achieved. The solution at the final time is shown in the right panel of Fig.~\ref{fig:graficoTest1}, while mass conservation and error norms are reported in Tab.~\ref{tab:errorTest3}.

\begin{table}[!h]
	\centering
	\begin{tabular}{|c|c|c|c|c|c|c|}
		\hline \multirow{2}{*}{$N_x$} & \multicolumn{3}{c|}{$\SL$} & \multicolumn{3}{c|}{$\SCOUT$} \\
		\cline{2-7} & $\Delta q$ & $\L^1$-err & $L^{1}$-ord & $\Delta q$ & $\L^{1}$-err & $\L^1$-ord\\ \hline \hline
		$50$ & $1.67e-17$ & $3.50e-03$ & - & $4.44e-18$ & $3.41e-03$ & - \\ \hline
		$100$ & $1.39e-18$ & $8.42e-04$ & $2.06$ & $5.25e-17$ & $8.42e-04$ & $2.02$ \\ \hline
		$200$ & $6.45e-17$ & $2.14e-04$ & $1.97$ & $1.10e-17$ & $2.10e-04$ & $2.01$ \\ \hline
		$400$ & $1.53e-16$ & $5.27e-05$ & $2.02$ & $1.09e-17$ & $5.24e-05$ & $2.00$ \\ \hline
		$800$ & $6.47e-17$ & $1.31e-05$ & $2.01$ & $8.53e-17$ & $1.31e-05$ & $2.00$ \\ \hline
	\end{tabular}
	\caption{Test 2: mass conservation, errors and accuracy. Comparison between $\SCOUT$ and $\SL$ schemes.}
	\label{tab:errorTest3}
\end{table}

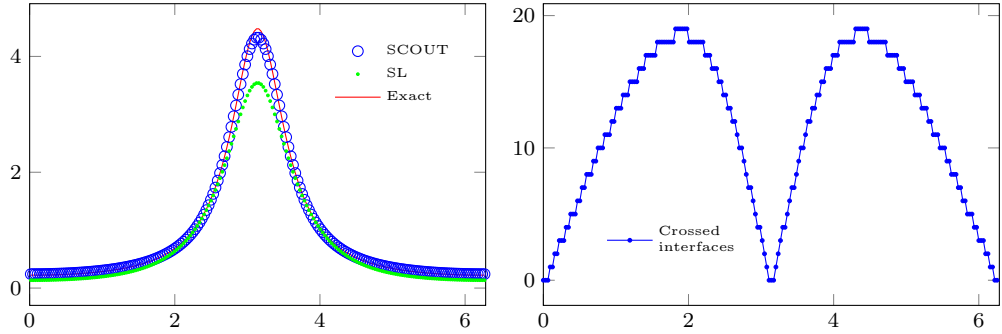
\begin{figure}[!h]
	\centering
	\begin{tikzpicture}

\pgfplotstableread{test4_convStep3.txt}{\graphData}

\begin{axis}[%
  width=6cm,height=4cm,scale only axis,
  xmin=0,xmax=2*pi,
  legend columns=1, 
  legend cell align={left},
  legend style={
    fill=none,
  	at={(0.95,0.9)},
  	draw=none,
    font=\tiny,
  }
]

\addplot[blue,mark=o,only marks] table[x=xcoor,y=solSLc] {\graphData};
\addlegendentry{$\SCOUT$};

\addplot[color=green,mark=*,only marks,mark size=0.5pt] table[x=xcoor,y=solSL] {\graphData};
\addlegendentry{$\SL$};

\addplot[color=red,mark=none] table[x=xcoor,y=solExa] {\graphData};
\addlegendentry{Exact};

\end{axis}
\end{tikzpicture}
	\begin{tikzpicture}

\pgfplotstableread{test4_convStep3_crossed.txt}{\graphData}

\begin{axis}[%
  width=6cm,height=4cm,scale only axis,
  xmin=0,xmax=2*pi,
  legend columns=1, 
  legend cell align={left},
  legend style={
    cells={align=left},
    fill=none,
  	at={(0.45,0.3)},
  	draw=none,
    font=\tiny,
  }
]

\addplot[blue,mark=*,mark size=0.7pt] table[x=xcoorInt,y=nCrossed] {\graphData};
\addlegendentry{Crossed\\interfaces};

\end{axis}
\end{tikzpicture}
	\caption{Variable coefficient case. Left: comparison among the $\SCOUT$ scheme, the $\SL$ scheme and the exact solution. Right: plot of the number of cell interfaces crossed by each trajectory during the first time step of the computation. Results are obtained by setting $N_x=200$, $T=1.5$, $\nu=0$, $\CFL=20$.}
	\label{fig:graficoTest4}
\end{figure}

\paragraph{Test 3: Variable coefficient case.} We continue by considering a linear advection equation with variable coefficient, that is 
\begin{equation}\label{eq:test4}
    \pder{t}{q} + \pder{x}{(\sin(x)q)} = 0,
\end{equation}
on the domain $X=[0,2\pi]$, setting $q_0(x)=1$, $T=1.5$, $\CFL=20$, and for which the exact solution is given by
\begin{equation}
    q(x,t) = \frac{\sin(2\arctan(e^{-t}\tan(\frac{x}{2})))}{\sin(x)}.
\end{equation}
The solution of \eqref{eq:test4} is advected with velocity $u(x)=\sin(x)$, causing the characteristics to be curves, rather than straight lines, and thus introducing a temporal error in the approximation of their foot. For this reason, we expect the classical $\SL$ scheme to be first order accurate, while also non-conservative, whereas the novel $\SCOUT$ scheme ensures second order of accuracy and mass conservation, as stated in Tab.~\ref{tab:errorTest4}. The final numerical solution is depicted in the left panel of Fig.~\ref{fig:graficoTest4}. We point out that, in this example, the high value of the $\CFL$ number causes the characteristics to cross a very large number of computational cells. As an example, the right panel of Fig.~\ref{fig:graficoTest4} shows, for each interface $\ipm$, $i=0,...,N_x$, the number of cell interfaces crossed by the Lagrangian trajectory in the first time step.

\begin{table}[!h]
	\centering
	\begin{tabular}{|c|c|c|c|c|c|c|}
		\hline \multirow{2}{*}{$N_x$} & \multicolumn{3}{c|}{$\SL$} & \multicolumn{3}{c|}{$\SCOUT$} \\
		\cline{2-7} & $\Delta q$ & $\L^1$-err & $L^{1}$-ord & $\Delta q$ & $\L^{1}$-err & $\L^1$-ord\\ \hline \hline
		$50$ & $2.49e+00$ & $2.83e+00$ & - & $0.00e+00$ & $1.60e+00$ & - \\ \hline
		$100$ & $1.99e+00$ & $2.14e+00$ & $0.40$ & $0.00e+00$ & $9.43e-01$ & $0.76$ \\ \hline
		$200$ & $8.97e-01$ & $8.97e-01$ & $1.26$ & $8.88e-16$ & $1.25e-01$ & $2.92$ \\ \hline
		$400$ & $4.30e-01$ & $4.30e-01$ & $1.06$ & $8.88e-16$ & $2.25e-02$ & $2.47$ \\ \hline
		$800$ & $2.12e-01$ & $2.12e-01$ & $1.02$ & $1.78e-15$ & $5.37e-03$ & $2.07$ \\ \hline
	\end{tabular}
	\caption{Test 3: mass conservation, errors and accuracy. While the $\SCOUT$ scheme achieves second order of accuracy and conserves the mass, the classical $\SL$ one is shown to be only first order accurate and non-conservative. $\CFL=20$.}
	\label{tab:errorTest4}
\end{table}

\begin{figure}[!h]
	\centering
	\begin{tikzpicture}

\pgfplotstableread{test9_convStep2.txt}{\graphData}

\begin{axis}[%
  width=6cm,height=4cm,scale only axis,
  xmin=-5,xmax=5,
  legend columns=1, 
  legend cell align={left},
  legend style={
    fill=none,
  	at={(0.5,0.9)},
  	draw=none,
    font=\tiny,
  }
]

\addplot[blue,mark=o,only marks] table[x=xcoor,y=solSLc] {\graphData};
\addlegendentry{$\SCOUT$};

\addplot[color=green,mark=*,only marks,mark size=0.5pt] table[x=xcoor,y=solSL] {\graphData};
\addlegendentry{$\SL$};

\addplot[color=red,mark=none] table[x=xcoor,y=solExa] {\graphData};
\addlegendentry{Exact};

\end{axis}
\end{tikzpicture}
	\begin{tikzpicture}

\pgfplotstableread{test9_convStep2_crossed.txt}{\graphData}

\begin{axis}[%
  width=6cm,height=4cm,scale only axis,
  xmin=-5,xmax=5,
  legend columns=1, 
  legend cell align={left},
  legend style={
    cells={align=left},
    fill=none,
  	at={(0.9,0.9)},
  	draw=none,
    font=\tiny,
  }
]

\addplot[blue,mark=*,mark size=0.7pt] table[x=xcoorInt,y=nCrossed] {\graphData};
\addlegendentry{Crossed\\interfaces};

\end{axis}
\end{tikzpicture}
	\caption{Burgers' equation with analytical solution: Tavelli test. Left: comparison among the $\SCOUT$ scheme, the $\SL$ scheme and the exact solution.Right: plot of the number of cell interfaces crossed by each trajectory during the first time step of the computation. Results are obtained by setting $N_x=100$, $T=0.2$, $\nu=0$, $\CFL=10$.}
	\label{fig:graficoTest9}
\end{figure}
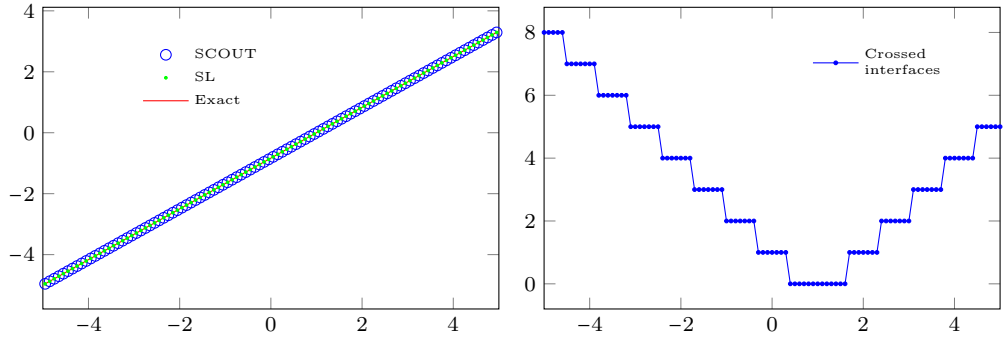

\paragraph{Test 4: Burgers' equation with linear data} In the sequel, we deal with the nonlinear inviscid Burgers' equation
\begin{equation}\label{eq:burgers}
    q_t + \left(\frac{q^2}{2}\right)_x = 0,
\end{equation}
whose solution is advected with velocity $u(q)=q$. As a starting point, we consider the following initial data:
\begin{equation}\label{eq:tavellisInitialData}
    q_0(x,t) = \alpha(t) + \beta x,
\end{equation}
where
\begin{equation}\label{eq:tavellisAlpha}
    \alpha(t) = \frac{\delta}{1+\beta t},
\end{equation}
being $\beta=1$ and $\delta=-1$ two parameters of the problem. For \eqref{eq:burgers}--\eqref{eq:tavellisInitialData}, an exact solution is given by
\begin{equation}\label{eq:tavellisExact}
    q(x,t) = \alpha(t) + \frac{\beta x}{1+\beta t},
\end{equation}
thus we can equip \eqref{eq:burgers} with Dirichlet boundary conditions. We verify the effectiveness of our method in the computational domain $X=[-5,5]$, with $T=0.2$ and $\CFL=10$. Similarly to the previous test, we show in Fig.~\ref{fig:graficoTest9} the final solution and the number of cells crossed by the characteristics during the first time step. Since the solution is linear in the spatial coordinate $x$, and the characteristics are straight lines, both schemes, $\SCOUT$ and $\SL$, approximate the solution at the machine precision error, as expected.

\paragraph{Test 5: Burgers' equation with sinusoidal initial data.} Another example for Burgers' equation \eqref{eq:burgers} is the one proposed in \cite{ConsSL2021} that assumes the following initial condition
\begin{equation}
    q_0(x) = \frac{\sqrt{2}}{2} + \sin(\pi x)
\end{equation}
and periodic boundary conditions in the domain $X=[0,2]$ and $\CFL=10$. Next, we also consider $\CFL=100$. The solution is advanced until different final times, in order to distinguish the case before and after the shock formation.

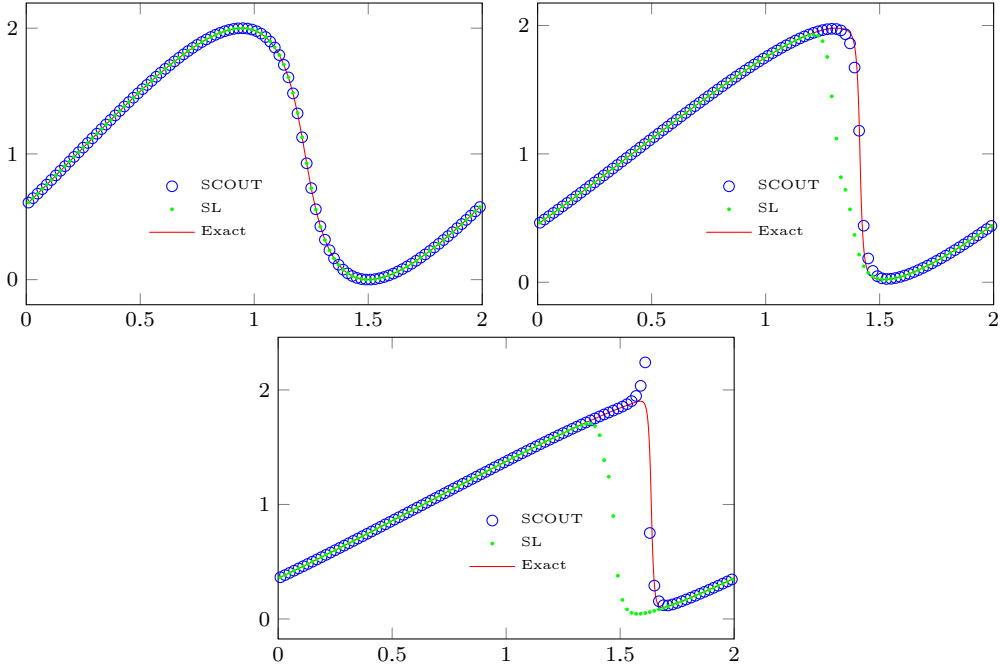
\begin{figure}[!h]
	\centering
	\begin{tikzpicture}

\pgfplotstableread{test6_convStep2_tFin1.txt}{\graphData}

\begin{axis}[%
  width=6cm,height=4cm,scale only axis,
  xmin=0,xmax=2,
  legend columns=1, 
  legend cell align={left},
  legend style={
    fill=none,
  	at={(0.55,0.45)},
  	draw=none,
    font=\tiny,
  }
]

\addplot[blue,mark=o,only marks] table[x=xcoor,y=solSLc] {\graphData};
\addlegendentry{$\SCOUT$};

\addplot[color=green,mark=*,only marks,mark size=0.5pt] table[x=xcoor,y=solSL] {\graphData};
\addlegendentry{$\SL$};

\addplot[color=red,mark=none] table[x=xcoor,y=solExa] {\graphData};
\addlegendentry{Exact};

\end{axis}
\end{tikzpicture}
	\begin{tikzpicture}

\pgfplotstableread{test6_convStep2_tFin2.txt}{\graphData}

\begin{axis}[%
  width=6cm,height=4cm,scale only axis,
  xmin=0,xmax=2,
  legend columns=1, 
  legend cell align={left},
  legend style={
    fill=none,
  	at={(0.65,0.45)},
  	draw=none,
    font=\tiny,
  }
]

\addplot[blue,mark=o,only marks] table[x=xcoor,y=solSLc] {\graphData};
\addlegendentry{$\SCOUT$};

\addplot[color=green,mark=*,only marks,mark size=0.5pt] table[x=xcoor,y=solSL] {\graphData};
\addlegendentry{$\SL$};

\pgfplotstableread{test6_convStep5_refSol_tFin2.txt}{\graphData}

\addplot[color=red,mark=none] table[x=xcoorRef,y=solRef] {\graphData};
\addlegendentry{Exact};

\end{axis}
\end{tikzpicture}
	\begin{tikzpicture}

\pgfplotstableread{test6_convStep2_tFin3.txt}{\graphData}

\begin{axis}[%
  width=6cm,height=4cm,scale only axis,
  xmin=0,xmax=2,
  legend columns=1, 
  legend cell align={left},
  legend style={
    fill=none,
  	at={(0.7,0.45)},
  	draw=none,
    font=\tiny,
  }
]

\addplot[blue,mark=o,only marks] table[x=xcoor,y=solSLc] {\graphData};
\addlegendentry{$\SCOUT$};

\addplot[color=green,mark=*,only marks,mark size=0.5pt] table[x=xcoor,y=solSL] {\graphData};
\addlegendentry{$\SL$};

\pgfplotstableread{test6_convStep5_refSol_tFin3.txt}{\graphData}

\addplot[color=red,mark=none] table[x=xcoorRef,y=solRef] {\graphData};
\addlegendentry{Exact};

\end{axis}
\end{tikzpicture}
	\caption{Numerical solutions for the nonlinear Burgers' equation with initial data $q_0(x)=\sqrt{2}/2+\sin(\pi x)$ computed on $N_x=100$ cells, with $\CFL=10$. Top left: $T=0.7/\pi$. Top right: $T=1.3/\pi$. Bottom: $T=2/\pi$. In the first case the solution is still smooth and no artificial viscosity is added. In the other two cases, we make a comparison with a reference solution computed on $N_x=800$ cells, considering also an additional artificial viscosity with coefficient $\nu=\num{5e-3}$.}
	\label{fig:graficoTest6}
\end{figure}

We first consider $T=0.7/\pi$, that is when the shock has not been formed yet, and we check the properties of the proposed scheme. As shown in Tab.~\ref{tab:errorTest6}, the $\SCOUT$ scheme results to be second order accurate and conservative, while the classical $\SL$ method achieves first order of accuracy, and is not conservative. We then increase the final time in order to see the shock formation and access the performance of the novel scheme in terms of mass conservation. Fig.~\ref{fig:graficoTest6} shows the results with $T=0.7\pi$, $T=1.3\pi$ and $T=2\pi$, where in the last two cases and additional artificial viscosity with constant coefficient $\nu=\num{5e-3}$ is added, and a reference solution with $N_x=800$ is computed. One can notice that the shock is well resolved by the conservative scheme, while it is completely misplaced by the classical $\SL$ method.

\begin{table}[!h]
	\centering
	\begin{tabular}{|c|c|c|c|c|c|c|}
		\hline \multirow{2}{*}{$N_x$} & \multicolumn{3}{c|}{$\SL$} & \multicolumn{3}{c|}{$\SCOUT$} \\
		\cline{2-7} & $\Delta q$ & $\L^1$-err & $L^{1}$-ord & $\Delta q$ & $\L^{1}$-err & $\L^1$-ord\\ \hline \hline
		$50$ & $2.35e-04$ & $3.17e-03$ & - & $2.22e-16$ & $1.98e-03$ & - \\ \hline
		$100$ & $1.50e-05$ & $1.02e-03$ & $1.63$ & $2.22e-16$ & $4.09e-04$ & $2.27$ \\ \hline
		$200$ & $1.68e-05$ & $3.83e-04$ & $1.42$ & $4.44e-16$ & $9.47e-05$ & $2.11$ \\ \hline
		$400$ & $5.54e-07$ & $1.45e-04$ & $1.40$ & $8.88e-16$ & $2.32e-05$ & $2.03$ \\ \hline
		$800$ & $3.87e-08$ & $6.78e-05$ & $1.10$ & $6.66e-16$ & $5.82e-06$ & $1.99$ \\ \hline
	\end{tabular}
	\caption{Burgers' equation with sinusoidal initial data: mass conservation, errors and accuracy. $\CFL=10$ and $T=0.7/\pi$. Comparison between $\SCOUT$ and $\SL$ schemes.}
	\label{tab:errorTest6}
\end{table}

\begin{figure}[!h]
	\centering
	\begin{tikzpicture}

\pgfplotstableread{test6_convStep3_CFL100_tFin1.txt}{\graphData}

\begin{axis}[%
  width=6cm,height=4cm,scale only axis,
  xmin=0,xmax=2,
  legend columns=1, 
  legend cell align={left},
  legend style={
    fill=none,
  	at={(0.55,0.45)},
  	draw=none,
    font=\tiny,
  }
]

\addplot[blue,mark=o,only marks] table[x=xcoor,y=solSLc] {\graphData};
\addlegendentry{$\SCOUT$};

\addplot[color=green,mark=*,only marks,mark size=0.5pt] table[x=xcoor,y=solSL] {\graphData};
\addlegendentry{$\SL$};

\addplot[color=red,mark=none] table[x=xcoor,y=solExa] {\graphData};
\addlegendentry{Exact};

\end{axis}
\end{tikzpicture}
	\begin{tikzpicture}

\pgfplotstableread{test6_convStep3_CFL100_tFin3.txt}{\graphData}

\begin{axis}[%
  width=6cm,height=4cm,scale only axis,
  xmin=0,xmax=2,
  legend columns=1, 
  legend cell align={left},
  legend style={
    fill=none,
  	at={(0.5,0.9)},
  	draw=none,
    font=\tiny,
  }
]

\addplot[blue,mark=o,only marks] table[x=xcoor,y=solSLc] {\graphData};
\addlegendentry{$\SCOUT$};

\addplot[color=green,mark=*,only marks,mark size=0.5pt] table[x=xcoor,y=solSL] {\graphData};
\addlegendentry{$\SL$};

\pgfplotstableread{test6_convStep5_refSol_CFL100_tFin3.txt}{\graphData}

\addplot[color=red,mark=none] table[x=xcoorRef,y=solRef] {\graphData};
\addlegendentry{Exact};

\end{axis}
\end{tikzpicture}
	\caption{Numerical solutions for the nonlinear Burgers' equation with initial data $q_0(x)=\sqrt{2}/2+\sin(\pi x)$ computed on $N_x=200$ cells, with $\CFL=100$. Left: $T=0.7/\pi$. Right: $T=2/\pi$. A reference solution computed on $N_x=800$ cells is also depicted, and an additional artificial viscosity with viscosity coefficient $\nu=\num{5e-3}$ is added in the second case.}
	\label{fig:graficoTest6_CFL100}
\end{figure}
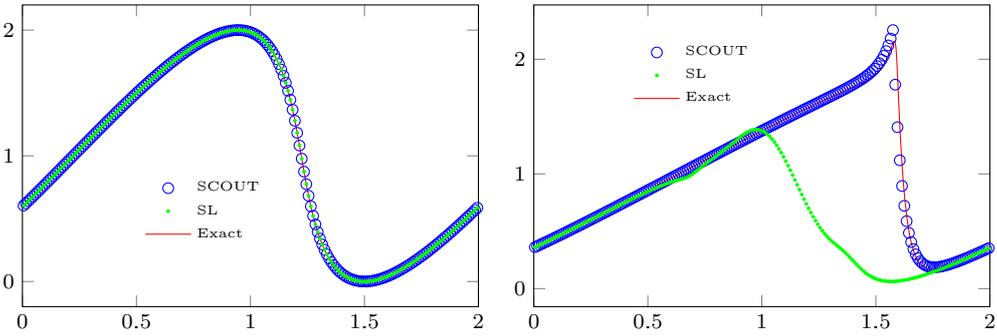

We are interested in testing our method with truly large values of the $\CFL$ number. In Fig.~\ref{fig:graficoTest6_CFL100} we show the results obtained by setting $\CFL=100$ at times $T=0.7\pi$ and $T=2\pi$. Moreover, Tab.~\ref{tab:errorTest6_CFL100} confirms the conservation property as well as the accuracy of the scheme, even in this extreme case.

\begin{table}[!h]
	\centering
	\begin{tabular}{|c|c|c|c|c|c|c|}
		\hline \multirow{2}{*}{$N_x$} & \multicolumn{3}{c|}{$\SL$} & \multicolumn{3}{c|}{$\SCOUT$} \\
		\cline{2-7} & $\Delta q$ & $\L^1$-err & $L^{1}$-ord & $\Delta q$ & $\L^{1}$-err & $\L^1$-ord\\ \hline \hline
		$50$ & $7.14e-05$ & $2.42e-03$ & - & $4.44e-16$ & $1.84e-03$ & - \\ \hline
		$100$ & $2.80e-06$ & $6.15e-04$ & $1.98$ & $0.00e+00$ & $4.05e-04$ & $2.19$ \\ \hline
		$200$ & $8.27e-06$ & $1.63e-04$ & $1.92$ & $0.00e+00$ & $9.31e-05$ & $2.12$ \\ \hline
		$400$ & $2.36e-07$ & $3.84e-05$ & $2.08$ & $4.44e-16$ & $2.24e-05$ & $2.05$ \\ \hline
		$800$ & $3.20e-08$ & $1.31e-05$ & $1.55$ & $1.33e-15$ & $5.55e-06$ & $2.02$ \\ \hline
		$1600$ & $8.74e-08$ & $5.35e-06$ & $1.29$ & $0.00e+00$ & $1.40e-06$ & $1.98$ \\ \hline
		$3200$ & $1.26e-09$ & $2.32e-06$ & $1.21$ & $8.88e-16$ & $3.59e-07$ & $1.96$ \\ \hline
	\end{tabular}
	\caption{Burgers' equation with sinusoidal initial data: mass conservation, errors and accuracy. $\CFL=100$ and $T=0.7/\pi$. Comparison between $\SCOUT$ and $\SL$ schemes.}
	\label{tab:errorTest6_CFL100}
\end{table}

\begin{figure}[!h]
	\centering
	\begin{tikzpicture}

\pgfplotstableread{test8_convStep3_CFL100_tFin1.txt}{\graphData}

\begin{axis}[%
  width=6cm,height=4cm,scale only axis,
  xmin=-5,xmax=5,
  legend columns=1, 
  legend cell align={left},
  legend style={
    fill=none,
  	at={(0.5,0.9)},
  	draw=none,
    font=\tiny,
  }
]

\addplot[blue,mark=o,only marks] table[x=xcoor,y=solSLc] {\graphData};
\addlegendentry{$\SCOUT$};

\addplot[color=green,mark=*,only marks,mark size=0.5pt] table[x=xcoor,y=solSL] {\graphData};
\addlegendentry{$\SL$};

\addplot[color=red,mark=none] table[x=xcoor,y=solExa] {\graphData};
\addlegendentry{Exact};

\end{axis}
\end{tikzpicture}
	\begin{tikzpicture}

\pgfplotstableread{test8_convStep3_CFL100_tFin2.txt}{\graphData}

\begin{axis}[%
  width=6cm,height=4cm,scale only axis,
  xmin=-5,xmax=5,
  legend columns=1, 
  legend cell align={left},
  legend style={
    fill=none,
  	at={(0.5,0.9)},
  	draw=none,
    font=\tiny,
  }
]

\addplot[blue,mark=o,only marks] table[x=xcoor,y=solSLc] {\graphData};
\addlegendentry{$\SCOUT$};

\addplot[color=green,mark=*,only marks,mark size=0.5pt] table[x=xcoor,y=solSL] {\graphData};
\addlegendentry{$\SL$};

\pgfplotstableread{test8_convStep5_refSol_CFL100_tFin2.txt}{\graphData}

\addplot[color=red,mark=none] table[x=xcoorRef,y=solRef] {\graphData};
\addlegendentry{Exact};

\end{axis}
\end{tikzpicture}
	\caption{Numerical solutions for the nonlinear Burgers' equation with Gaussian data $q_0(x)=1 + 0.5 e^{-2 x^2}$ with $N_x=200$ and $\CFL=100$. Left: $T=0.9$. Right: $T=2$. In the first case the solution is still smooth and no artificial viscosity is added. In the second case, we make a comparison against a reference solution computed on $N_x=800$ cells, considering also an additional artificial viscosity with viscosity coefficient $\nu=\num{5e-3}$.}
	\label{fig:graficoTest8_CFL100}
\end{figure}
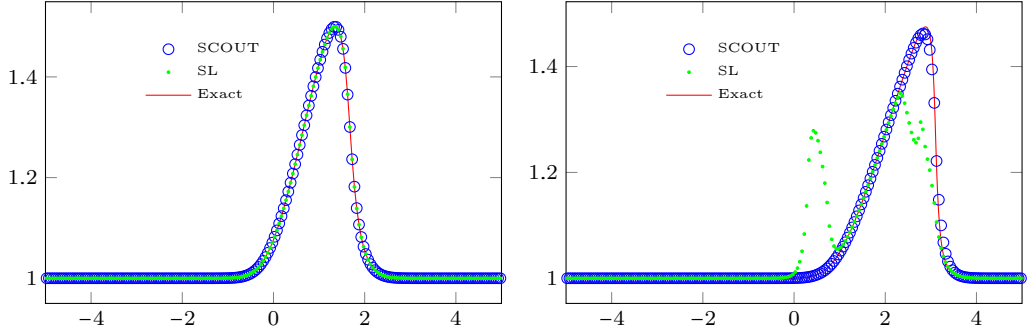

\paragraph{Test 6: Burgers' equation with Gaussian initial data.} Similarly to the previous test, we now consider as initial data the Gaussian profile
\begin{equation}
    q_0(x) = 1 + \frac{1}{2}e^{-2 x^2}
\end{equation}
in the domain $X=[-5,5]$, with periodic boundary conditions, at different final times $T=0.9$ and $T=2$, before and after the shock formation, respectively. We run the simulation with $\CFL=100$. The mass conservation and accuracy results are shown in Tab.~\ref{tab:errorTest8_CFL100}, while plots of the solution are depicted in Fig.~\ref{fig:graficoTest8_CFL100}, confirming once again both accuracy and mass conservation.  

\begin{table}[!h]
	\centering
	\begin{tabular}{|c|c|c|c|c|c|c|}
		\hline \multirow{2}{*}{$N_x$} & \multicolumn{3}{c|}{$\SL$} & \multicolumn{3}{c|}{$\SCOUT$} \\
		\cline{2-7} & $\Delta q$ & $\L^1$-err & $L^{1}$-ord & $\Delta q$ & $\L^{1}$-err & $\L^1$-ord\\ \hline \hline
		$50$ & $5.66e-04$ & $1.23e-02$ & - & $1.78e-15$ & $8.43e-03$ & - \\ \hline
		$100$ & $2.76e-04$ & $2.67e-03$ & $2.20$ & $0.00e+00$ & $1.71e-03$ & $2.30$ \\ \hline
		$200$ & $2.15e-05$ & $6.41e-04$ & $2.06$ & $0.00e+00$ & $3.78e-04$ & $2.18$ \\ \hline
		$400$ & $1.69e-06$ & $1.49e-04$ & $2.11$ & $0.00e+00$ & $8.90e-05$ & $2.09$ \\ \hline
		$800$ & $6.23e-06$ & $5.82e-05$ & $1.35$ & $0.00e+00$ & $2.16e-05$ & $2.04$ \\ \hline
		$1600$ & $6.06e-07$ & $1.78e-05$ & $1.71$ & $8.88e-15$ & $5.34e-06$ & $2.01$ \\ \hline
		$3200$ & $2.88e-07$ & $6.58e-06$ & $1.43$ & $0.00e+00$ & $1.34e-06$ & $2.00$ \\ \hline
	\end{tabular}
	\caption{Burgers' equation with Gaussian initial profile: mass conservation, errors and accuracy. $\CFL=100$ and $T=0.9$. Comparison between $\SCOUT$ and $\SL$ schemes.}
	\label{tab:errorTest8_CFL100}
\end{table}

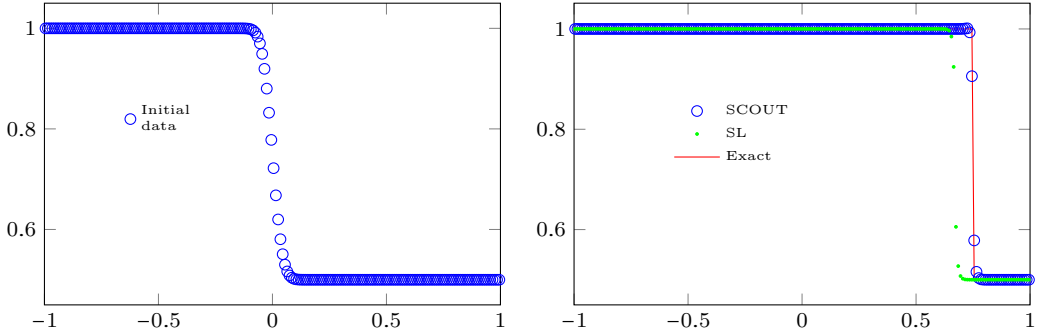
\begin{figure}[!h]
	\centering
	\begin{tikzpicture}

\pgfplotstableread{test7_convStep3_iniData.txt}{\graphData}

\begin{axis}[%
  width=6cm,height=4cm,scale only axis,
  xmin=-1,xmax=1,
  legend columns=1, 
  legend cell align={left},
  legend style={
    cells={align=left},
    fill=none,
  	at={(0.35,0.7)},
  	draw=none,
    font=\tiny,
  }
]

\addplot[blue,mark=o,only marks] table[x=xcoorIni,y=solIni] {\graphData};
\addlegendentry{Initial\\data};

\end{axis}
\end{tikzpicture}
	\begin{tikzpicture}

\pgfplotstableread{test7_convStep3.txt}{\graphData}

\begin{axis}[%
  width=6cm,height=4cm,scale only axis,
  xmin=-1,xmax=1,
  legend columns=1, 
  legend cell align={left},
  legend style={
    fill=none,
  	at={(0.5,0.7)},
  	draw=none,
    font=\tiny,
  }
]

\addplot[blue,mark=o,only marks] table[x=xcoor,y=solSLc] {\graphData};
\addlegendentry{$\SCOUT$};

\addplot[color=green,mark=*,only marks,mark size=0.5pt] table[x=xcoor,y=solSL] {\graphData};
\addlegendentry{$\SL$};

\addplot[color=red,mark=none] table[x=xcoor,y=solExa] {\graphData};
\addlegendentry{Exact};

\end{axis}
\end{tikzpicture}
	\caption{Test 7: $N_x=200$, $T=1$, $\nu=\num{5e-4}$, $\CFL=10$. The initial profile is shown on the left, while on the right we plot the solution at the final time for both the $\SCOUT$ and the $\SL$ scheme as well as the exact one.}
	\label{fig:graficoTest7}
\end{figure}

\paragraph{Test 7: Traveling shock wave.} We conclude the test suite for the advective case by considering the well-known traveling shock wave that arises when solving the inviscid Burgers' equation \eqref{eq:burgers}, setting as initial profile
\begin{equation}
    q_0(x) = 0.5 \left( q_L+q_R + (q_R-q_L)\,\erf\left(\frac{x}{0.05}\right) \right),
\end{equation}
where $q_L=1$, $q_R=0.5$ and $\erf()$ is the error function. Despite the regularity of the initial data, the solution becomes discontinuous, causing the conservation property to be crucial for the accuracy of the solution and the correct shock capturing property in terms of both magnitude and position. Fig.~\ref{fig:graficoTest7} shows the initial data and the traveling shock wave at time $T=1$. The exact position of the shock is known, thus we can access the performance of the $\SCOUT$ scheme, with respect to the classical $\SL$ one.

\begin{figure}[!h]
	\centering
	\begin{tikzpicture}

\pgfplotstableread{test11_convStep4_nu02_CNvisc.txt}{\graphData}

\begin{axis}[%
  width=6cm,height=4cm,scale only axis,
  xmin=-3,xmax=3,
  legend columns=1, 
  legend cell align={left},
  legend style={
    fill=none,
  	at={(0.5,0.4)},
  	draw=none,
    font=\tiny,
  }
]

\addplot[blue,mark=o,only marks] table[x=xcoor,y=solSLc] {\graphData};
\addlegendentry{$\SCOUTSLCN$};

\addplot[color=green,mark=*,only marks,mark size=0.5pt] table[x=xcoor,y=solSL] {\graphData};
\addlegendentry{$\SL$};

\addplot[color=red,mark=none] table[x=xcoor,y=solExa] {\graphData};
\addlegendentry{exact};

\end{axis}
\end{tikzpicture}
	\begin{tikzpicture}

\pgfplotstableread{test11_convStep4_nu02_CNviscClassic.txt}{\graphData}

\begin{axis}[%
  width=6cm,height=4cm,scale only axis,
  xmin=-3,xmax=3,
  legend columns=1, 
  legend cell align={left},
  legend style={
    fill=none,
  	at={(0.5,0.4)},
  	draw=none,
    font=\tiny,
  }
]

\addplot[blue,mark=o,only marks] table[x=xcoor,y=solSLc] {\graphData};
\addlegendentry{$\SCOUTCN$};

\addplot[color=green,mark=*,only marks,mark size=0.5pt] table[x=xcoor,y=solSL] {\graphData};
\addlegendentry{$\SL$};

\addplot[color=red,mark=none] table[x=xcoor,y=solExa] {\graphData};
\addlegendentry{exact};

\end{axis}
\end{tikzpicture}
	\caption{Test 8: $N_x=400$, $T=1$, $\nu=\num{1e-2}$, $\CFL=10$. Comparison among the final solution obtained with the $\SCOUT$ and $\SL$ method, and the exact solution. Left: $\SCOUT$ scheme with semi-Lagrangian Crank-Nicolson method for the diffusion terms. Right: $\SCOUT$ scheme with classical Crank-Nicolson method for the diffusion terms.}
	\label{fig:graficoTest11}
\end{figure}
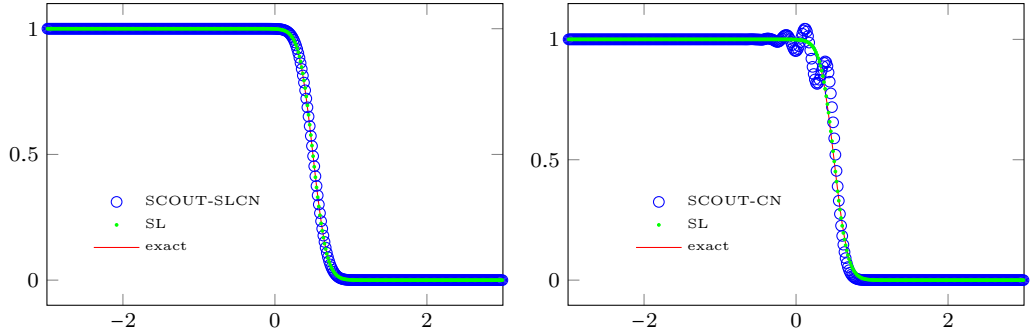

\paragraph{Test 8: Advection-diffusion with constant velocity.} In the last three test cases, we want to benchmark the novel $\SCOUT$ method for the model problem \eqref{eq:advdiff1d}, as it has been designed in Sec.~\ref{sec:advdiff}. We start by considering the advection-diffusion process occurring at constant velocity $\Bar{u}=1$, for $\nu=\num{1e-2}$, and choosing as initial data a step function that becomes gradually smoothed during the evolution due to the diffusion process. We set $X=[-3,3]$, $T=1$, and
\begin{equation}
    q_0(x) =
    \begin{cases}
        1 \quad \text{ for } x<-1/2\\
        0 \quad \text{ for } x\geq-1/2
    \end{cases}.
\end{equation}
An exact solution for this problem is given by
\begin{equation}\label{eq:exaAdvDiffErf}
    q(x,t) = \frac{1}{2}-\frac{1}{2}\erf\left( \frac{x+1/2-t}{\sqrt{4\nu t}} \right),
\end{equation}
when $t>0$. The initial data is actually set from \eqref{eq:exaAdvDiffErf}, choosing $t=0.1$ as initial profile to avoid loss of regularity in the initial condition. To highlight the fundamental difference between a classical decoupled Crank-Nicolson scheme and the new semi-Lagrangian Crank-Nicolson scheme for the diffusion terms detailed in Sec.~\ref{sec:advdiff} and embedded in the $\SCOUT$ method, this test case is run in both configurations. We use the names $\SCOUTSLCN$ and $\SCOUTCN$ to refer to these two schemes, respectively. From Tab.~\ref{tab:errorTest11_CNvisc} and Tab.~\ref{tab:errorTest11_CNviscClassic} one can notice that the $\SCOUT$ method achieves second order of accuracy, while the classical Crank-Nicolson scheme does not. This is due to the fact that in this latter case, for each computational node $x_i$, information are taken only at $(x_i,\tn)$ and $(x_i,\tnn)$, thus affecting the accuracy and the quality of the solution, as it is clearly depicted in Fig.~\ref{fig:graficoTest11}.

\begin{table}[!h]
	\centering
	\begin{tabular}{|c|c|c|c|c|}
		\hline \multirow{2}{*}{$N_x$} & \multicolumn{2}{c|}{$\SL$} & \multicolumn{2}{c|}{$\SCOUTSLCN$} \\
		\cline{2-5} & $\L^1$-err & $L^{1}$-ord & $\L^{1}$-err & $\L^1$-ord\\ \hline \hline
		$50$ & $2.94e-02$ & - & $9.12e-03$ & - \\ \hline
		$100$ & $1.26e-02$ & $1.22$ & $1.77e-03$ & $2.36$ \\ \hline
		$200$ & $6.49e-03$ & $0.96$ & $4.28e-04$ & $2.05$ \\ \hline
		$400$ & $3.34e-03$ & $0.96$ & $1.05e-04$ & $2.02$ \\ \hline
		$800$ & $1.70e-03$ & $0.98$ & $2.70e-05$ & $1.96$ \\ \hline
		$1600$ & $8.59e-04$ & $0.99$ & $6.80e-06$ & $1.99$ \\ \hline
		$3200$ & $4.32e-04$ & $0.99$ & $1.72e-06$ & $1.99$ \\ \hline
	\end{tabular}
	\caption{Advection-diffusion equation with constant velocity solved with the $\SCOUTSLCN$ method: errors and accuracy. $\CFL=10$, $\nu=\num{1e-2}$, and $T=1$. Comparison between $\SCOUTSLCN$ and $\SL$ schemes.}
	\label{tab:errorTest11_CNvisc}
\end{table}

\begin{table}[!h]
	\centering
	\begin{tabular}{|c|c|c|c|c|c|c|}
		\hline \multirow{2}{*}{$N_x$} & \multicolumn{2}{c|}{$\SL$} & \multicolumn{2}{c|}{$\SCOUTCN$} \\
		\cline{2-5} & $\L^1$-err & $L^{1}$-ord & $\L^{1}$-err & $\L^1$-ord\\ \hline \hline
		$50$ & $2.94e-02$ & - & $8.63e-02$ & - \\ \hline
		$100$ & $1.26e-02$ & $1.22$ & $7.08e-02$ & $0.29$ \\ \hline
		$200$ & $6.13e-03$ & $1.04$ & $7.28e-02$ & $-0.04$ \\ \hline
		$400$ & $3.41e-03$ & $0.85$ & $5.79e-02$ & $0.33$ \\ \hline
		$800$ & $1.70e-03$ & $1.01$ & $1.94e-02$ & $1.58$ \\ \hline
		$1600$ & $8.59e-04$ & $0.98$ & $8.44e-03$ & $1.20$ \\ \hline
		$3200$ & $4.32e-04$ & $0.99$ & $4.16e-03$ & $1.02$ \\ \hline
	\end{tabular}
	\caption{Advection-diffusion equation with constant velocity solved with the classical Crank-Nicolson scheme for the diffusion terms: errors and accuracy. $\CFL=10$, $\nu=\num{1e-2}$, and $T=1$. Comparison between $\SCOUTCN$ and $\SL$ schemes.}
	\label{tab:errorTest11_CNviscClassic}
\end{table}

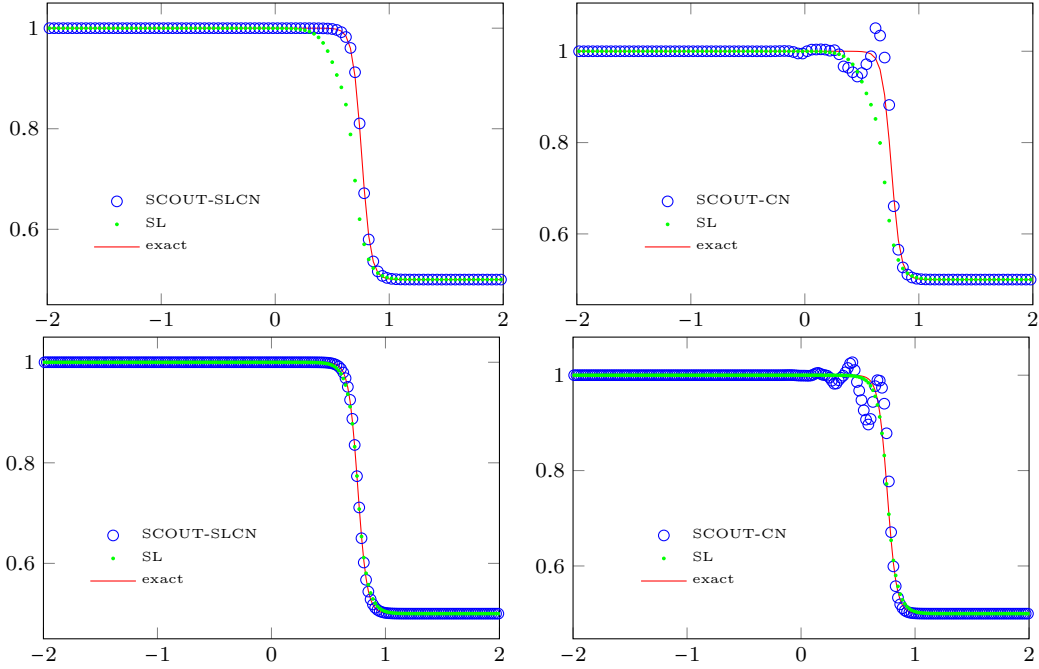
\begin{figure}[!h]
	\centering
	\begin{tikzpicture}

\pgfplotstableread{test14_convStep2_CNvisc.txt}{\graphData}

\begin{axis}[%
  width=6cm,height=4cm,scale only axis,
  xmin=-2,xmax=2,
  legend columns=1, 
  legend cell align={left},
  legend style={
    fill=none,
  	at={(0.5,0.4)},
  	draw=none,
    font=\tiny,
  }
]

\addplot[blue,mark=o,only marks] table[x=xcoor,y=solSLc] {\graphData};
\addlegendentry{$\SCOUTSLCN$};

\addplot[color=green,mark=*,only marks,mark size=0.5pt] table[x=xcoor,y=solSL] {\graphData};
\addlegendentry{$\SL$};

\addplot[color=red,mark=none] table[x=xcoor,y=solExa] {\graphData};
\addlegendentry{exact};

\end{axis}
\end{tikzpicture}
	\begin{tikzpicture}

\pgfplotstableread{test14_convStep2_CNviscClassic.txt}{\graphData}

\begin{axis}[%
  width=6cm,height=4cm,scale only axis,
  xmin=-2,xmax=2,
  legend columns=1, 
  legend cell align={left},
  legend style={
    fill=none,
  	at={(0.5,0.4)},
  	draw=none,
    font=\tiny,
  }
]

\addplot[blue,mark=o,only marks] table[x=xcoor,y=solSLc] {\graphData};
\addlegendentry{$\SCOUTCN$};

\addplot[color=green,mark=*,only marks,mark size=0.5pt] table[x=xcoor,y=solSL] {\graphData};
\addlegendentry{$\SL$};

\addplot[color=red,mark=none] table[x=xcoor,y=solExa] {\graphData};
\addlegendentry{exact};

\end{axis}
\end{tikzpicture}
	\begin{tikzpicture}

\pgfplotstableread{test14_convStep3_CNvisc.txt}{\graphData}

\begin{axis}[%
  width=6cm,height=4cm,scale only axis,
  xmin=-2,xmax=2,
  legend columns=1, 
  legend cell align={left},
  legend style={
    fill=none,
  	at={(0.5,0.4)},
  	draw=none,
    font=\tiny,
  }
]

\addplot[blue,mark=o,only marks] table[x=xcoor,y=solSLc] {\graphData};
\addlegendentry{$\SCOUTSLCN$};

\addplot[color=green,mark=*,only marks,mark size=0.5pt] table[x=xcoor,y=solSL] {\graphData};
\addlegendentry{$\SL$};

\addplot[color=red,mark=none] table[x=xcoor,y=solExa] {\graphData};
\addlegendentry{exact};

\end{axis}
\end{tikzpicture}
	\begin{tikzpicture}

\pgfplotstableread{test14_convStep3_CNviscClassic.txt}{\graphData}

\begin{axis}[%
  width=6cm,height=4cm,scale only axis,
  xmin=-2,xmax=2,
  legend columns=1, 
  legend cell align={left},
  legend style={
    fill=none,
  	at={(0.5,0.4)},
  	draw=none,
    font=\tiny,
  }
]

\addplot[blue,mark=o,only marks] table[x=xcoor,y=solSLc] {\graphData};
\addlegendentry{$\SCOUTCN$};

\addplot[color=green,mark=*,only marks,mark size=0.5pt] table[x=xcoor,y=solSL] {\graphData};
\addlegendentry{$\SL$};

\addplot[color=red,mark=none] table[x=xcoor,y=solExa] {\graphData};
\addlegendentry{exact};

\end{axis}
\end{tikzpicture}
	\caption{Test 9: Top row: $N_x=100$, $T=1$, $\nu=0.01$, $\CFL=10$. Comparison among the final solution obtained with the $\SCOUT$ and $\SL$ method, and the exact solution. Left: $\SCOUT$ scheme with semi-Lagrangian Crank-Nicolson method for diffusion terms. Right: $\SCOUT$ scheme with classical Crank-Nicolson method for the diffusion terms. Bottom row: same configurations using $N_x=200$. Note that the mass is not conserved by the $\SL$ scheme. High values of the $\CFL$ number produce undesirable behaviors in the decoupled case on the right.}
	\label{fig:graficoTest14}
\end{figure}

\paragraph{Test 9: Travelling wave solution for viscous Burgers' equation.} Moving to the nonlinear case, we consider again the viscous Burgers' equation
\begin{equation}\label{eq:viscBurgers}
    q_t + \left(\frac{q^2}{2}\right)_x = \nu q_{xx},
\end{equation}
in the domain $X=[-2,2]$, coupled with the initial condition
\begin{equation}
    q_0(x) =
    \begin{cases}
        q_L \quad \text{ for } x<0\\
        q_R \quad \text{ for } x\geq0
    \end{cases},
\end{equation}
with $q_L=1$ and $q_R=0.5$. Similarly as done previously, $q_0$ is smoothed out by \eqref{eq:viscBurgers}, so that the exact solution takes the form $q^{\nu}(x,t) = w(x-st)$, where $s=(q_L+q_R)/2$ and
\begin{equation}
    w(y) = q_R + \frac{1}{2}(q_L-q_R)\left[ 1-\tanh\left( \frac{q_L-q_R}{4\nu} y \right) \right].
\end{equation}
This type of problem is very relevant when dealing with conservation laws, especially regarding the uniqueness of their solution, and we refer the reader to \cite{leveque1992} for a comprehensive explanation. We set $\nu=\num{1e-2}$, $\CFL=10$ and let the solution evolve up to $T=1$, starting from $t_0=0.01$. Tab.~\ref{tab:errorTest14_CNvisc} and Tab.~\ref{tab:errorTest14_CNviscClassic} are coherent with our considerations about the semi-Lagrangian discretization of the diffusion term, as for the linear case. From Fig.~\ref{fig:graficoTest14} it is also noticeable that the profile is not well-captured by the classic $SL$ scheme, due to the lack of mass conservation.

\begin{table}[!h]
	\centering
	\begin{tabular}{|c|c|c|c|c|}
		\hline \multirow{2}{*}{$N_x$} & \multicolumn{2}{c|}{$\SL$} & \multicolumn{2}{c|}{$\SCOUTSLCN$} \\
		\cline{2-5} & $\L^1$-err & $L^{1}$-ord & $\L^{1}$-err & $\L^1$-ord\\ \hline \hline
		$50$ & $1.31e-01$ & - & $1.67e-02$ & - \\ \hline
		$100$ & $5.07e-02$ & $1.37$ & $2.17e-03$ & $2.94$ \\ \hline
		$200$ & $3.98e-03$ & $3.67$ & $7.70e-04$ & $1.49$ \\ \hline
		$400$ & $1.92e-03$ & $1.05$ & $1.68e-04$ & $2.20$ \\ \hline
		$800$ & $9.44e-04$ & $1.02$ & $4.01e-05$ & $2.06$ \\ \hline
		$1600$ & $4.68e-04$ & $1.01$ & $9.85e-06$ & $2.03$ \\ \hline
		$3200$ & $2.33e-04$ & $1.01$ & $2.44e-06$ & $2.01$ \\ \hline
	\end{tabular}
	\caption{Viscous Burgers' equation with travelling wave solution: errors and accuracy solved with the $\SCOUTSLCN$ method. $\CFL=10$, $\nu=\num{1e-2}$, and $T=1$. Comparison between $\SCOUTSLCN$ and $\SL$ schemes.}
	\label{tab:errorTest14_CNvisc}
\end{table}

\begin{table}[!h]
	\centering
	\begin{tabular}{|c|c|c|c|c|c|c|}
		\hline \multirow{2}{*}{$N_x$} & \multicolumn{2}{c|}{$\SL$} & \multicolumn{2}{c|}{$\SCOUTCN$} \\
		\cline{2-5} & $\L^1$-err & $L^{1}$-ord & $\L^{1}$-err & $\L^1$-ord\\ \hline \hline
		$50$ & $1.31e-01$ & - & $4.50e-02$ & - \\ \hline
		$100$ & $4.86e-02$ & $1.43$ & $2.65e-02$ & $0.77$ \\ \hline
		$200$ & $3.81e-03$ & $3.67$ & $2.41e-02$ & $0.13$ \\ \hline
		$400$ & $1.87e-03$ & $1.03$ & $1.39e-02$ & $0.80$ \\ \hline
		$800$ & $9.43e-04$ & $0.99$ & $5.28e-03$ & $1.39$ \\ \hline
		$1600$ & $4.68e-04$ & $1.01$ & $2.57e-03$ & $1.04$ \\ \hline
		$3200$ & $2.33e-04$ & $1.01$ & $1.29e-03$ & $1.00$ \\ \hline
	\end{tabular}
	\caption{Viscous Burgers' equation with travelling wave solution: errors and accuracy solved with the classical Crank-Nicolson scheme for the diffusion terms. $\CFL=10$, $\nu=\num{1e-2}$, and $T=1$. Comparison between $\SCOUTCN$ and $\SL$ schemes.}
	\label{tab:errorTest14_CNviscClassic}
\end{table}

\begin{figure}[!h]
	\centering
	\begin{tikzpicture}

\pgfplotstableread{test10_convStep2.txt}{\graphData}

\begin{axis}[%
  width=6cm,height=4cm,scale only axis,
  xmin=0,xmax=2*pi,
  legend columns=1, 
  legend cell align={left},
  legend style={
    fill=none,
  	at={(0.6,0.4)},
  	draw=none,
    font=\tiny,
  }
]

\addplot[blue,mark=o,only marks] table[x=xcoor,y=solSLc] {\graphData};
\addlegendentry{$\SCOUTSLCN$};

\addplot[color=green,mark=*,only marks,mark size=0.5pt] table[x=xcoor,y=solSL] {\graphData};
\addlegendentry{$\SL$};

\addplot[color=red,mark=none] table[x=xcoor,y=solExa] {\graphData};
\addlegendentry{exact};

\end{axis}
\end{tikzpicture}
	\begin{tikzpicture}

\pgfplotstableread{test10_convStep3.txt}{\graphData}

\begin{axis}[%
  width=6cm,height=4cm,scale only axis,
  xmin=0,xmax=2*pi,
  legend columns=1, 
  legend cell align={left},
  legend style={
    fill=none,
  	at={(0.6,0.4)},
  	draw=none,
    font=\tiny,
  }
]

\addplot[blue,mark=o,only marks] table[x=xcoor,y=solSLc] {\graphData};
\addlegendentry{$\SCOUTSLCN$};

\addplot[color=green,mark=*,only marks,mark size=0.5pt] table[x=xcoor,y=solSL] {\graphData};
\addlegendentry{$\SL$};

\addplot[color=red,mark=none] table[x=xcoor,y=solExa] {\graphData};
\addlegendentry{exact};

\end{axis}
\end{tikzpicture}
	\caption{Test 8: final solution of the viscous Burgers' equation at time $T=1/\pi$ and $\CFL=10$, with $N_x=100$ (left) and $N_x=200$ (right). The conservation of mass is guaranteed by the $\SCOUT$ scheme and not by the classical $\SL$ method. In both cases, the semi-Lagrangian approach is embedded in the Crank-Nicolson scheme.}
	\label{fig:graficoTest10}
\end{figure}
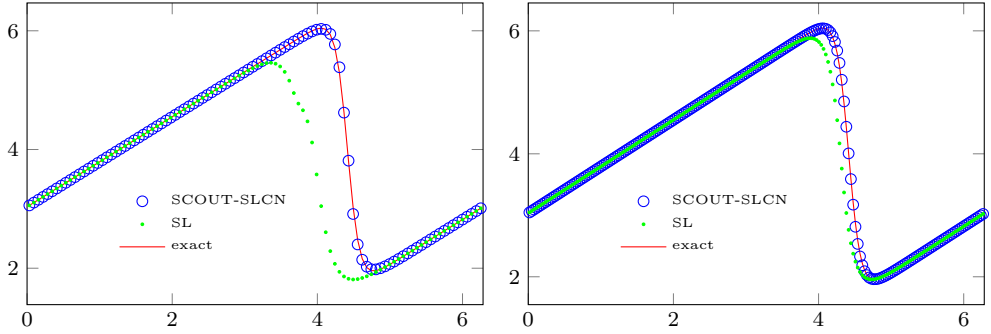

\paragraph{Test 10: Viscous Burgers' equation with periodic data.} To conclude, we consider again the model problem \eqref{eq:viscBurgers} with periodic boundary conditions, in the domain $X=[0,2\pi]$, setting as initial condition
\begin{equation}
    q_0(x) = 4-2\nu\frac{\pder{x}{\phi(x,0)}}{\phi(x,0)},
\end{equation}
where
\begin{equation}
    \phi(x,t) = e^{-\frac{x-4t}{4\nu(t+1)}} + e^{-\frac{x-4t-2\pi}{4\nu(t+1)}},
\end{equation}
and for which an accurate reference solution can be computed following \cite{advdiffExa1986}. We set $\nu=\num{2e-1}$, $T=1/\pi$, $\CFL=10$ and observe the results obtained with the $\SCOUTSLCN$ method, compared to the classical $\SL$. In Fig.~\ref{fig:graficoTest10} the solution at the final time are depicted for $N_x=100$ and $N_x=200$, from left to right, respectively, and one can immediately appreciate that the conservation of mass, as well as the correct profile, is only guaranteed by the novel $\SCOUTSLCN$ method. This result is quantitatively confirmed in Tab.~\ref{tab:errorTest10}, together with the accuracy orders. 

\begin{table}[!h]
	\centering
	\begin{tabular}{|c|c|c|c|c|c|c|}
		\hline \multirow{2}{*}{$N_x$} & \multicolumn{3}{c|}{$\SL$} & \multicolumn{3}{c|}{$\SCOUTSLCN$} \\
		\cline{2-7} & $\Delta q$ & $\L^1$-err & $L^{1}$-ord & $\Delta q$ & $\L^{1}$-err & $\L^1$-ord\\ \hline \hline
		$50$ & $4.15e+00$ & $4.15e+00$ & - & $3.55e-15$ & $6.17e-01$ & - \\ \hline
		$100$ & $2.21e+00$ & $2.21e+00$ & $0.91$ & $6.04e-14$ & $3.70e-02$ & $4.06$ \\ \hline
		$200$ & $5.11e-01$ & $5.19e-01$ & $2.09$ & $1.42e-14$ & $2.02e-02$ & $0.87$ \\ \hline
		$400$ & $6.08e-06$ & $5.31e-02$ & $3.29$ & $4.87e-13$ & $4.60e-03$ & $2.13$ \\ \hline
		$800$ & $5.50e-05$ & $2.62e-02$ & $1.02$ & $1.63e-12$ & $1.03e-03$ & $2.16$ \\ \hline
		$1600$ & $3.98e-05$ & $1.30e-02$ & $1.01$ & $1.00e-11$ & $2.49e-04$ & $2.05$ \\ \hline
		$3200$ & $2.14e-05$ & $6.49e-03$ & $1.00$ & $4.29e-11$ & $6.16e-05$ & $2.01$ \\ \hline
	\end{tabular}
	\caption{Viscous Burgers' equation with periodic data: mass conservation, errors and accuracy. $\CFL=10$ and $T=1/\pi$. Comparison between $\SCOUTSLCN$ and $\SL$ schemes.}
	\label{tab:errorTest10}
\end{table}

\section{Conclusions}\label{sec:concl}
In this paper, we have designed, implemented, and validated the $\SCOUT$ method for advection-diffusion problems. The scheme relies on a conservative semi-Lagrangian discretization of both advection and diffusion, achieves second-order accuracy in space and time, and is unconditionally stable. 
{The advective contribution is evaluated upon integration of the governing PDE over a space-time domain bounded by a characteristic curve and two segments aligned with the spatial and the time axis: subsequent application of Gauss' theorem eventually enables to compute the time integrals of the flux only using spatial information at the beginning of the timestep. In the nonlinear case only, a Newton method is employed to approximate the foot of the characteristic and the integral along the characteristic might not vanish, leading to an additional term to be computed. In all cases, mass conservation in ensured by construction.}
Diffusion is incorporated within the conservative semi-Lagrangian framework through a characteristic-based Crank-Nicolson time-marching scheme, that preserves second-order accuracy and unconditional stability, which is not the case if a classical Crank-Nicolson method, not based on the semi-Lagrangian formulation nor coupled with the advection, is adopted. Several linear and nonlinear advection and advection-diffusion test cases are presented. We have systematically compared the new $\SCOUT$ method against the classical semi-Lagrangian scheme solving the governing PDE in non-conservative form in terms of accuracy and robustness. $\CFL$ numbers up to 100 have been used to deeply test the novel scheme.

Future developments include extending the proposed numerical method to nonlinear systems of conservation laws, such as the Euler-Fourier model, and to two-dimensional Cartesian and unstructured meshes via multidimensional space-time integration on the control volumes to guarantee conservation. We also plan to investigate higher-order space–time conservative semi-Lagrangian schemes based on the $\SCOUT$ framework.

\section*{Acknowledgments}
WB received financial support by the Italian Ministry of University and Research (MUR) with the PRIN Project 2022 No. 2022N9BM3N.
MS received financial support by the Italian Ministry of University and Research (MUR) with the PRIN Project 2022 No. 2022JH87B4.
The authors are members of the GNCS-INdAM (\textit{Istituto Nazionale di Alta Matematica}) group.

	\section*{Declarations}

	\paragraph{Funding.} WB received financial support by the Italian Ministry of University (MUR) with the PRIN Project 2022 No. 2022N9BM3N. MS received financial support by the Italian Ministry of University and Research (MUR) with the PRIN Project 2022 No. 2022JH87B4. The authors are also grateful to the SHARK-FV 2025 workshop, where the main ideas of this paper were first developed.

	\paragraph{Conflicts of interest.} The authors declare that they have no conflict of interest.

	\paragraph{Availability of data and material.} Data and material are available upon reasonable request addressed to the corresponding author.

	\paragraph{Code availability.} The code is written in \texttt{Matlab} programming language and is available upon reasonable request addressed to the corresponding author.

	\paragraph{Ethics approval.} Not applicable.

	\paragraph{Consent to participate.} Not applicable.

	\paragraph{Consent for publication.} Not applicable.

	\section*{Data availability}
	The datasets generated during the current study are available from the corresponding author upon reasonable request.

	\bibliographystyle{plain}
	\bibliography{biblio}   %

\end{document}